\documentclass[a4paper]{article}%
\usepackage{amssymb}
\usepackage{amsfonts}
\usepackage{amsmath}
\usepackage[a4paper]{geometry}
\usepackage{wrapfig}
\usepackage{graphicx}%
\providecommand{\U}[1]{\protect\rule{.1in}{.1in}}
\newtheorem{theorem}{Theorem}

\begin{document}

\title{Universal Hyperbolic Geometry II: A pictorial overview}
\author{N J Wildberger\\School of Mathematics and Statistics\\UNSW Sydney 2052 Australia}
\date{}
\maketitle

\begin{abstract}
This article provides a simple pictorial introduction to universal hyperbolic
geometry. We explain how to understand the subject using only elementary
projective geometry, augmented by a distinguished circle. This provides a
completely algebraic framework for hyperbolic geometry, valid over the
rational numbers (and indeed any field not of characteristic two), and gives
us many new and beautiful theorems. These results are accurately illustrated
with colour diagrams, and the reader is invited to check them with ruler
constructions and measurements.

\end{abstract}

\section{Introduction}

This paper introduces \textit{hyperbolic geometry} using only elementary
mathematics, without any analysis, and in particular without transcendental
functions. Classical hyperbolic geometry, (see for example \cite{Greenberg},
\cite{Hartshorne}, \cite{Iversen}, \cite{Lenz}, \cite{McCleary},
\cite{Milnor}, \cite{Stillwell}, \cite{Ungar}, \cite{Wolfe}), is usually an
advanced topic studied in the senior years of a university mathematics
program, often built up from a foundation of differential geometry. In recent
years, a new, simpler and completely algebraic understanding of this subject
has emerged, building on the ideas of \textit{rational trigonometry
}(\cite{Wild1} and \cite{Wild2}). This approach is called \textit{universal
hyperbolic geometry}, because it extends the theory to more general settings,
namely to arbitrary fields (usually avoiding characteristic two), and because
it generalizes to other quadratic forms (see \cite{Wild4}).

The basic reference is \textit{Universal Hyperbolic Geometry I: Trigonometry}
(\cite{Wild5}), which contains accurate definitions, many formulas and
complete proofs, but no diagrams. This paper complements that one, providing a
pictorial introduction to the subject with a minimum of formulas and no
proofs, essentially relying only on planar projective geometry. The reader is
encouraged to verify theorems by making explicit constructions and
measurements; aside from a single base (null) circle, with only a ruler one
can check most of the assertions of this paper in special cases. Alternatively
a modern geometry program such as The Geometer's Sketchpad, C.A.R., Cabri,
GeoAlgebra or Cinderella illustrates the subject with a little effort.

Our approach extends the classical \textit{Cayley Beltrami Klein}, or
\textit{projective}, model of hyperbolic geometry, whose underlying space is
the interior of a disk, with lines being straight line segments. In our
formulation we consider also the boundary of the disk, which we call the
\textit{null circle}, also points outside the disk, and also points at
infinity. The lines are now complete lines in the sense of projective
geometry, not segments, and include also \textit{null lines} which are tangent
to the null circle, and lines which do not meet the null circle, including the
line at infinity. This orientation is familiar to classical geometers (see for
example \cite{Benz}, \cite{Brauner}, \cite{Coolidge}, \cite{Coxeter},
\cite{Sommerville}), but it is not well-known to students because of the
current dominance of the differential geometric point of view. A novel aspect
of this paper is that we introduce our metrical concepts---\textit{quadrance},
\textit{spread }and \textit{quadrea}---purely in projective terms. It means
that only high school algebra suffices to set up the subject and make
computations. The proofs however rely generally on computer calculations
involving polynomial or rational function identities; some may be found in
(\cite{Wild5}), others will appear elsewhere.

Many of the results are illustrated with two diagrams, one illustrating the
situation in the classical setting using interior points of the null circle,
and another with more general points. The fundamental metrical notions of
\textit{quadrance} between points and \textit{spread }between lines are
undefined when null points or null lines are involved, but most theorems
involving them apply equally to points and lines interior or exterior of the
null circle. The reader should be aware that in more advanced work the
distinction between these two types of points and lines also becomes
significant. Instead of area of hyperbolic triangles, we work with a rational
analog called \textit{quadrea}.

With this algebraic approach we can develop geometry over the \textit{rational
numbers}---in my view, always the most important field. The natural
connections between geometry and number theory are then not suppressed, but
enrich both subjects. Later in the series we will also illustrate hyperbolic
geometry over \textit{finite fields, }in the direction of \cite{Angel} and
\cite{Terras}, where \textit{counting} becomes important.

\section{The projective plane}

Hyperbolic geometry may be visualized as the geometry of the projective plane,
augmented by a distinguished circle $c$ (in fact a more general conic may also
be used). Since projective geometry is not these days as familiar as it was in
former times, we begin by reviewing some of the basic notions. The starting
point is the affine plane---familiar from Euclidean geometry and Cartesian
coordinate geometry---containing the usual points $\left[  x,y\right]  $ and
(straight) lines with equations $ax+by=c$. The affine plane is augmented by
introducing a \textit{new point }for every family of parallel lines. In this
introductory section we use\textit{\ parallel} in the usual sense of Euclidean
geometry, so that the lines with equations $a_{1}x+b_{1}y=c_{1}$ and
$a_{2}x+b_{2}y=c_{2}$ are parallel precisely when $a_{1}b_{2}-a_{2}b_{1}=0.$
Later on we will see that there is a different, hyperbolic meaning of
`parallel' (which is different from the usage in classical hyperbolic
geometry!) The new point, one for each family of parallel lines, is a
`\textit{point at infinity}'. We also introduce one \textit{new line}, the
`\textit{line at infinity'}, which passes through every point at infinity.

Algebraically the projective plane may be defined without reference to the
Euclidean plane, with points specified by homogeneous coordinates, or
proportions, of the form $\left[  x:y:z\right]  $. Points $\left[
x:y:1\right]  $ correspond to the affine plane, and points at infinity are of
the form $\left[  x:y:0\right]  $. The lines also are specified by homogeneous
coordinates, now of the form $\left(  a:b:c\right)  $, with the pairing
between the point $\left[  x:y:z\right]  $ and the line $\left(  a:b:c\right)
$ given by%
\begin{equation}
ax+by-cz=0. \label{LineEquationDef}%
\end{equation}

This particular relation is the characterizing equation for hyperbolic
geometry; for spherical/elliptic geometry a different convention between
points and lines is used, where the line $\left\langle a:b:c\right\rangle $
passes through the point $\left[  x:y:z\right]  $ precisely when $ax+by+cz=0.$
Note that we use round brackets for lines in hyperbolic geometry.

We will visualize the projective plane as an extension of the affine plane,
with the usual property that any two distinct points $a$ and $b$ determine
exactly one line which passes through them both, called the \textbf{join} of
$a$ and $b,$ and denoted by $ab,$ and with the \textit{new} property that any
two lines $L$ and $M$ determine exactly one point which lies on them both,
called the \textbf{meet }of $L$ and $K,$ and denoted by $LM.$

In projective geometry, the notion of parallel lines disappears, since now
\textit{any} two lines meet. Familiar measurements, such as distance and
angle, are also absent. It is really the \textit{geometry of the straightedge}.

\begin{figure}[h]%
\centering
\includegraphics[
height=2.1096in,
width=3.9269in
]%
{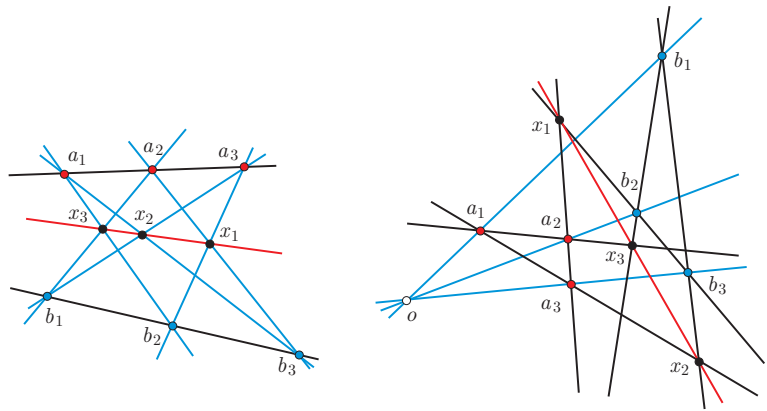}%
\caption{Theorems of Pappus and Desargues}%
\label{Pappus and Desargues}%
\end{figure}

Despite its historical importance, intrinsic beauty and simplicity, projective
geometry is these days sadly neglected in the school and university
curriculum. Perhaps the wider realization that it actually underpins
hyperbolic geometry will lead to a renaissance of the subject! Most readers
will know the two basic theorems in the subject, which are illustrated in
Figure \ref{Pappus and Desargues}.
%TCIMACRO{\FRAME{fhFU}{3.9269in}{2.1096in}{0pt}{\Qcb{Theorems of Pappus and
%Desargues}}{\Qlb{Pappus and Desargues}}{pappus_and_desargues.EPS}%
%width 3.9269in;  height 2.1096in;  depth 0pt;  original-width 5.1224in;
%original-height 2.7405in;  cropleft "0";  croptop "1";  cropright "1";
%cropbottom "0";  filename 'Pappus_and_Desargues.EPS';file-properties "XNPEU";}%
%} }%
%BeginExpansion
%EndExpansion

\textit{Pappus' theorem }asserts that if $a_{1}$, $a_{2}$ and $a_{3}$ are
collinear, and $b_{1},b_{2}$ and $b_{3}$ are collinear, then $x_{1}%
\equiv\left(  a_{2}b_{3}\right)  \left(  a_{3}b_{2}\right)  $, $x_{2}%
\equiv\left(  a_{3}b_{1}\right)  \left(  a_{1}b_{3}\right)  $ and $x_{3}%
\equiv\left(  a_{1}b_{2}\right)  \left(  a_{2}b_{1}\right)  $ are collinear.
\textit{Desargues' theorem }asserts that if $a_{1}b_{1}$, $a_{2}b_{2}$ and
$a_{3}b_{3}$ are concurrent, then $x_{1}\equiv\left(  a_{2}a_{3}\right)
\left(  b_{2}b_{3}\right)  $, $x_{2}\equiv\left(  a_{3}a_{1}\right)  \left(
b_{3}b_{1}\right)  $ and $x_{3}\equiv\left(  a_{1}a_{2}\right)  \left(
b_{1}b_{2}\right)  $ are collinear. This is often stated in the form that if
two triangles are perspective from a point, then they are also perspective
from a line.

A further important notion concerns four collinear points $a,b,c$ and $d$ on a
line $L,$ in any order. Suppose we choose affine coordinates on $L$ so that
the coordinates of $a,b,c$ and $d$ are respectively $x,y,z$ and $w.$ Then the
\textbf{cross-ratio} is defined to be the extended number (possibly $\infty$)
given by the ratio of ratios:
\[
\left(  a,b:c,d\right)  \equiv\left(  \frac{a-c}{b-c}\right)  /\left(
\frac{a-d}{b-d}\right)  .
\]
This is independent of the choice of affine coordinates on $L$.
%TCIMACRO{\FRAME{fhFU}{1.7335in}{1.3292in}{0pt}{\Qcb{Projective invariance of
%cross-ratio: $\left(  a,b:c,d\right)  =\left(  a_{1},b_{1}:c_{1},d_{1}\right)
%$}}{\Qlb{CrossRatio}}{crossratio.EPS}{\special{ language "Scientific Word";
%type "GRAPHIC";  maintain-aspect-ratio TRUE;  display "USEDEF";
%valid_file "F";  width 1.7335in;  height 1.3292in;  depth 0pt;
%original-width 2.4319in;  original-height 1.8568in;  cropleft "0";
%filename 'crossratio.EPS';file-properties "XNPEU";}} }%
%BeginExpansion
\begin{figure}[h]%
\centering
\includegraphics[
height=1.3292in,
width=1.7335in
]%
{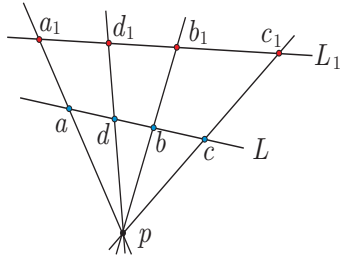}%
\caption{Projective invariance of cross-ratio: $\left(  a,b:c,d\right)
=\left(  a_{1},b_{1}:c_{1},d_{1}\right)  $}%
\label{CrossRatio}%
\end{figure}
%EndExpansion

Moreover it is also projectively invariant, meaning that if $a_{1},b_{1}%
,c_{1}$ and $d_{1}$ are also collinear points on a line $L_{1}$ which are
perspective to $a,b,c$ and $d$ from some point $p,$ as in Figure
\ref{CrossRatio}, then $\left(  a,b:c,d\right)  =\left(  a_{1},b_{1}%
:c_{1},d_{1}\right)  $.

The cross ratio is the most important invariant in projective geometry, and
will be the basis of \textit{quadrance} between points and \textit{spread
}between lines, but we give also analytic expressions for these quantities in
homogeneous coordinates; both are \textit{rational functions} of the inputs.
The actual values assumed by the cross-ratio depend ultimately on the field
over which our geometry is based, which is in principle quite arbitrary. To
start it helps to restrict our attention to the \textit{rational numbers},
invariably the most natural, familiar and important field. So we will adopt a
scientific approach, identifying our sheet of paper with (part of) the
\textit{rational number plane}.

Here are a few more basic definitions. A \textbf{side }$\overline{a_{1}a_{2}%
}=\left\{  a_{1},a_{2}\right\}  $ is a set of two points. A \textbf{vertex
}$\overline{L_{1}L_{2}}=\left\{  L_{1},L_{2}\right\}  $ is a set of two lines.
A \textbf{couple} $\overline{aL}=\left\{  a,L\right\}  $ is a set consisting
of a point $a$ and a line $L$. A \textbf{triangle }$\overline{a_{1}a_{2}a_{3}%
}=\left\{  a_{1},a_{2},a_{3}\right\}  $ is a set of three points which are not
collinear. A\textbf{\ trilateral }$\overline{L_{1}L_{2}L_{3}}=\left\{
L_{1},L_{2},L_{3}\right\}  $ is a set of three lines which are not concurrent.
Every triangle $\overline{a_{1}a_{2}a_{3}}$ has three sides, namely
$\overline{a_{1}a_{2}},$ $\overline{a_{2}a_{3}}$ and $\overline{a_{1}a_{3}}$,
and similarly any trilateral $\overline{L_{1}L_{2}L_{3}}$ has three vertices,
namely $\overline{L_{1}L_{2}},$ $\overline{L_{2}L_{3}}$ and $\overline
{L_{1}L_{3}}.$

Since the points of a triangle and the lines of a trilateral are distinct, any
triangle $\overline{a_{1}a_{2}a_{3}}$ determines an \textbf{associated
trilateral} $\overline{L_{1}L_{2}L_{3}}$ where $L_{1}\equiv a_{2}a_{3},$
$L_{2}\equiv a_{1}a_{3}$ and $L_{3}\equiv a_{1}a_{2}.$ Conversely any
trilateral $\overline{L_{1}L_{2}L_{3}}$ determines an \textbf{associated
triangle }$\overline{a_{1}a_{2}a_{3}}$ where $a_{1}\equiv L_{2}L_{3},$
$a_{2}\equiv L_{1}L_{3}$ and $a_{3}\equiv L_{1}L_{2}.$

\section{Duality via polarity}

We are now ready to introduce the \textbf{hyperbolic plane}, which is just the
projective plane which we have just been describing, augmented by a
distinguished Euclidean circle $c$ in this plane, called the \textbf{null
circle, }which appears in our diagrams always in blue. The points lying on $c$
have a distinguished role, and are called \textbf{null points}. The lines
tangent to $c$ have a distinguished role, and are called \textbf{null lines}.

All other points and lines, including the points at infinity and the line at
infinity, are for the purposes of elementary universal hyperbolic geometry
treated in a non-preferential manner. In particular we do \textit{not
}restrict our attention to only \textbf{interior points} lying inside the
circle $c$; this is a big difference with classical hyperbolic geometry;
\textbf{exterior points} lying outside the circle are equally important.
Similarly we do \textit{not} restrict our attention to only \textbf{interior
lines }which meet $c$ in two points; \textbf{exterior lines} which do not meet
$c$ are equally important. Note also that these notions can be defined purely
projectively once the null circle $c$ has been specified: interior points do
not lie on null lines, while exterior points do, and interior lines pass
through null points, while exterior lines do not.

For those who prefer to work with coordinates, we may choose our circle to
have homogeneous equation $x^{2}+y^{2}-z^{2}=0,$ or in the plane $z=1$ with
coordinates $X\equiv x/z$ and $Y\equiv y/z,$ simply the unit circle
$X^{2}+Y^{2}=1.$

The presence of the distinguished null circle $c$ has as its main consequence
a \textit{complete duality} between points and lines of the projective plane,
in the sense that every point $a$ has associated to it a particular line
$a^{\bot}$ and conversely. This duality is one of the ways in which universal
hyperbolic geometry is very different from classical hyperbolic geometry, and
it arises from a standard construction in projective geometry involving the
distinguished null circle $c$---the notion of \textit{polarity}. How polarity
defines duality is central to the subject.%
%TCIMACRO{\FRAME{fhFU}{10.5227cm}{4.9788cm}{0pt}{\Qcb{Duality and pole-polar
%pairs}}{\Qlb{PolesPolars}}{DualQuadBoth.eps}%
%width 10.5227cm;  height 4.9788cm;  depth 0pt;  original-width 13.8404cm;
%original-height 6.499cm;  cropleft "0";  croptop "1";  cropright "1";
%cropbottom "0";  filename 'DualQuadBoth.eps';file-properties "XNPEU";}} }%
%BeginExpansion
\begin{figure}[h]%
\centering
\includegraphics[
height=4.9788cm,
width=10.5227cm
]%
{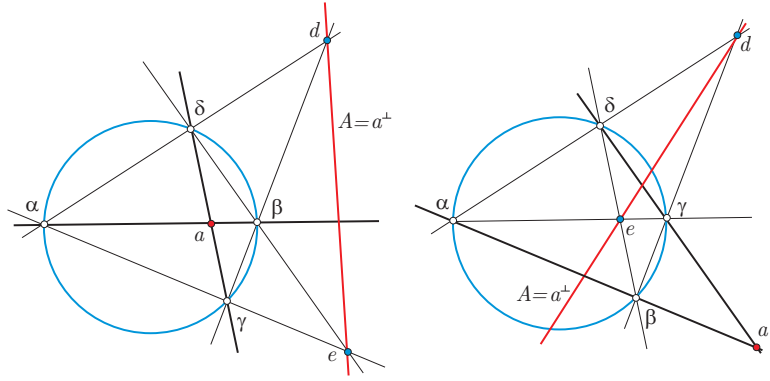}%
\caption{Duality and pole-polar pairs}%
\label{PolesPolars}%
\end{figure}
%EndExpansion

To define the polar of a point $a$ with respect to the null circle $c,$ draw
any two lines through $a$ which each both pass through $c$ in two distinct
points, say $\alpha,\beta$ and $\gamma,\delta$ respectively; this is always
possible. Now here is a beautiful fact from projective geometry:\textit{\ if
}$d$ and $e$ are the other two diagonal points of the quadrilateral
$\alpha\beta\gamma\delta,$ then \textit{the line }$de$ \textit{does not depend
on the two chosen lines through }$a,$\textit{\ but only on }$a$%
\textit{\ itself. }So we say that $de=A\equiv a^{\bot}$ is the \textbf{dual
line} of the point $a$, and conversely $a=A^{\bot}$ is the \textbf{dual point}
of the line $A.$

\begin{figure}[h]%
\centering
\includegraphics[
height=1.7974in,
width=2.0473in
]%
{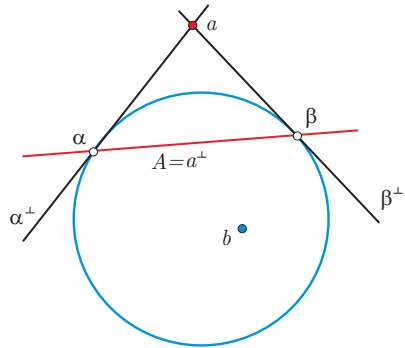}%
\caption{Null points $\alpha,\beta$ and their dual null lines}%
\label{NullPointLine}%
\end{figure}

The picture is as in Figure \ref{PolesPolars}, showing two different possible
configurations for which the above prescriptions both hold. In case $a$ is
external to the circle, it is also possible to construct $A\equiv a^{\bot}$
from the tangents to the null circle $c$ passing through $a$ as in Figure
\ref{NullPointLine}, but this does not work for an interior point such as $b.$

The construction shows that there is in fact a symmetry between the initial
point $a$ and the diagonal points $d$ and $e,$ so that if we started with the
point $d,$ its dual would be the line $ae,$ and if we started with the point
$e,$ its dual would be the point $ad.$ So this shows another fundamental fact:
\textit{if }$d$\textit{\ lies on the dual }$a^{\perp}$\textit{\ of the point
}$a,$\textit{\ then }$a$\textit{\ lies on the dual }$d^{\perp}$\textit{\ of
the point }$d.$\textit{\ }

So to invert the construction, given a line $A,$ choose two points $d$ and $e
$ on it, find the dual lines $d^{\perp}$ and $e^{\perp},$ and define $a\equiv
A^{\perp}=d^{\perp}e^{\perp}.$ It is at this point that we need the projective
plane, with its points and line at infinity, for if $d^{\perp}$ and $e^{\perp
}$ were Euclidean parallel, then $d^{\perp}e^{\perp}$ would be a point at
infinity. This situation occurs if we take our line $A$ to be a diameter of
the null circle $c,$ in the sense of Euclidean geometry.

What happens when the point $a$ is null? In that case the quadrilateral in the
above construction degenerates, and the polar line $a^{\perp}$ is the null
line tangent to the null circle at $a.$ So every null line is dual to the
unique null point which lies on it, as shown also in Figure
\ref{NullPointLine}.%
%TCIMACRO{\FRAME{fhFU}{2.0473in}{1.7974in}{0pt}{\Qcb{Null points $\alpha,\beta$
%and their dual null lines}}{\Qlb{NullPointLine}}{nullpointline.EPS}%
%width 2.0473in;  height 1.7974in;  depth 0pt;  original-width 2.7224in;
%original-height 2.386in;  cropleft "0";  croptop "1";  cropright "1";
%cropbottom "0";  filename 'NullPointLine.EPS';file-properties "XNPEU";}} }%
%BeginExpansion
%EndExpansion

The duality between points and lines is surprisingly simple to describe in
terms of homogeneous coordinates. The pole of the point $a\equiv\left[
x:y:z\right]  $ is just the line $a^{\perp}\equiv\left(  x:y:z\right)  $, so
duality amounts to simply changing from one kind of brackets to the other!
This explains why we chose the hyperbolic form of the pairing
(\ref{LineEquationDef}).

In the previous section we saw that every triangle $\overline{a_{1}a_{2}a_{3}%
}$ has an associated trilateral $\overline{L_{1}L_{2}L_{3}}$ and conversely.
Now we see that there is another natural trilateral associated to
$\overline{a_{1}a_{2}a_{3}},$ namely the \textbf{dual trilateral }%
$\overline{a_{1}^{\bot}a_{2}^{\bot}a_{3}^{\bot}}.$ Conversely to any
trilateral $\overline{L_{1}L_{2}L_{3}}$ there is associated the \textbf{dual
triangle} $\overline{L_{1}^{\bot}L_{2}^{\bot}L_{3}^{\bot}}$.
%TCIMACRO{\FRAME{fhFU}{2.9622in}{1.6239in}{0pt}{\Qcb{A triangle $\overline
%{a_{1}a_{2}a_{3}}$ and its dual trilateral $\overline{A_{1}A_{2}A_{3}}$}}%
%{}{dualtriangle.EPS}{\special{ language "Scientific Word";  type "GRAPHIC";
%width 2.9622in;  height 1.6239in;  depth 0pt;  original-width 5.7372in;
%original-height 3.1332in;  cropleft "0";  croptop "1";  cropright "1";
%cropbottom "0";  filename 'DualTriangle.EPS';file-properties "XNPEU";}} }%
%BeginExpansion
\begin{figure}[h]%
\centering
\includegraphics[
height=1.6239in,
width=2.9622in
]%
{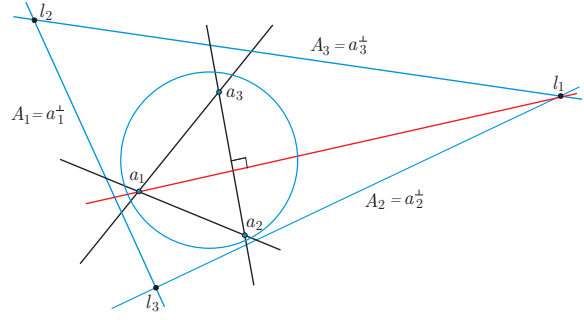}%
\caption{A triangle $\overline{a_{1}a_{2}a_{3}}$ and its dual trilateral
$\overline{A_{1}A_{2}A_{3}}$}%
\end{figure}
%EndExpansion

We also say that a couple $\overline{aL}$ is \textbf{dual} precisely when $a$
and $L$ are dual, and is otherwise \textbf{non-dual}.

So to summarize: \textit{in (planar) universal hyperbolic geometry, duality
implies that points and lines are treated completely symmetrically}. This is a
significant departure both from Euclidean geometry and from classical
hyperbolic geometry---in both of those theories, points and lines play quite
different roles. In universal hyperbolic geometry, the above duality principle
implies that every theorem can be dualized to create a (possibly) new theorem.
We will state many theorems together with their duals, but to keep the length
of this paper reasonable we will not extend this to all the theorems; the
reader is encouraged to find statements and draw pictures of the dual results
in these other cases.

\section{Perpendicularity}

The notions of perpendicular and parallel differ dramatically between
Euclidean and hyperbolic geometries. In the affine geometry on which Euclidean
geometry is based, the notion of \textit{parallel lines }is fundamental, while
\textit{perpendicular lines }are determined by a quadratic form. In hyperbolic
geometry, the situation is reversed---perpendicularity is more fundamental,
and in fact parallelism is defined in terms of it!

Another novel feature is that perpendicularity applies not only to lines, but
\textit{also to points}. This is a consequence of the fundamental duality we
have already established between points and lines.%
%TCIMACRO{\FRAME{fhFU}{4.0349in}{1.814in}{0pt}{\Qcb{Perpendicular points and
%lines}}{\Qlb{PerpendicularPointsLines}}{perpendicularlinesboth.eps}%
%width 4.0349in;  height 1.814in;  depth 0pt;  original-width 5.1872in;
%original-height 2.3177in;  cropleft "0";  croptop "1";  cropright "1";
%filename 'PerpendicularLinesBoth.eps';file-properties "XNPEU";}} }%
%BeginExpansion
\begin{figure}[h]%
\centering
\includegraphics[
height=1.814in,
width=4.0349in
]%
{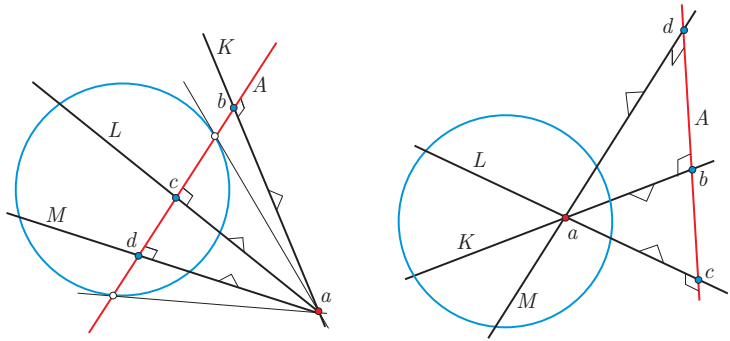}%
\caption{Perpendicular points and lines}%
\label{PerpendicularPointsLines}%
\end{figure}
%EndExpansion

Perpendicularity in the hyperbolic setting is easy to introduce once we have
duality. We say that a point $b$ is \textbf{perpendicular }to a point $a$
precisely when $b$ lies on the dual line $a^{\perp}$ of $a.$ This is
equivalent to $a$ lying on the dual line $b^{\perp}$ of $b,$ so the relation
is symmetric, and we write%
\[
a\perp b.
\]
Dually a line $L$ is perpendicular to a line $M$ precisely when $L$ passes
through the dual point $M^{\perp}$ of $M.$ This is equivalent to $M$ passing
through the dual point $L^{\perp}$ of $L,$ and we write%
\[
L\perp M.
\]

Figure \ref{PerpendicularPointsLines} shows our pictorial conventions for
perpendicularity: the line $A\ $is dual to the point $a$, so points $b,c$ and
$d$ lying on $A$ are perpendicular to $a,$ and this is recorded by a small
(right) corner placed on the join of the perpendicular points, and between
them. Also the lines $K,L$ and $M$ pass through $a,$ so are perpendicular to
$A,$ and this is recorded as usual by a small parallelogram at the meet of the
perpendicular lines.

Our first theorem records two basic facts that are obvious from the
definitions so far.

\begin{theorem}
[Altitude line and point]For any non-dual couple $\overline{aL}$, there is a
unique line $N$ which passes through $a$ and is perpendicular to $L,$ namely
$N\equiv aL^{\perp},$ and there is a unique point $n$ which lies on $L$ and is
perpendicular to $a, $ namely $n\equiv a^{\perp}L.$ Furthermore $N$ and $n$
are dual.%
%TCIMACRO{\FRAME{fhFU}{3.3031in}{1.4316in}{0pt}{\Qcb{Altitude line $N$ and
%altitude point $n$ of a couple $\overline{aL}$}}{\Qlb{AltitudeLinePoint}%
%}{altitudelinepoint2.EPS}{\special{ language "Scientific Word";
%type "GRAPHIC";  maintain-aspect-ratio TRUE;  display "USEDEF";
%valid_file "F";  width 3.3031in;  height 1.4316in;  depth 0pt;
%original-width 5.9205in;  original-height 2.5486in;  cropleft "0";
%filename 'AltitudeLinePoint2.EPS';file-properties "XNPEU";}} }%
%BeginExpansion
\begin{figure}[h]%
\centering
\includegraphics[
height=1.4316in,
width=3.3031in
]%
{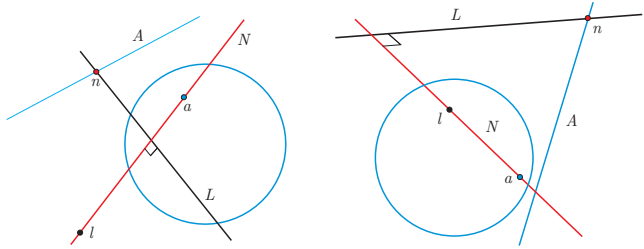}%
\caption{Altitude line $N$ and altitude point $n$ of a couple $\overline{aL}$}%
\label{AltitudeLinePoint}%
\end{figure}
%EndExpansion

\end{theorem}

We call $N$ the\textbf{\ altitude line} to $L$ through $a,$ and $n$ the
\textbf{altitude point }to $a$ on $L.$ In case $a$ and $L$ are dual, any line
through $a$ is perpendicular to $L,$ and any point lying on $L$ is
perpendicular to $a.$ While the former idea is familiar, the latter is not.
Figure \ref{AltitudeLinePoint} shows that if we restrict ourselves to the
inside of the null circle, altitude points are invisible, so it is no surprise
that the concept is missing from classical hyperbolic geometry. Remember that
we are obliged to respect the balance which duality provides us!

If a triangle $\overline{a_{1}a_{2}a_{3}}$ has say $a_{1}$ dual to $a_{2}%
a_{3},$ then any line through $a_{1}$ will be perpendicular to the opposite
line $a_{2}a_{3},$ and we say the triangle is \textbf{dual}. A triangle is
\textbf{non-dual }precisely when each of its points is not dual to the
opposite line. Similar definitions apply to trilaterals.

Somewhat surprisingly, the next result is \textit{not true} in classical
hyperbolic geometry---a very conspicuous absence that too often goes
unmentioned in books on the subject! The reason is that the orthocenter of a
triangle of interior points might well be an exterior point, as the left
diagram in Figure \ref{OrthocenterBoth} shows. The absence of a distinguished
orthocenter partially explains why the study of triangles in classical
hyperbolic geometry is relatively undeveloped. With universal hyperbolic
geometry, triangle geometry enters a rich new phase.

\begin{theorem}
[Orthocenter and ortholine]The altitude lines of a non-dual triangle meet at a
unique point $o$, called the \textbf{orthocenter} of the triangle. The
altitude points of a non-dual trilateral join along a unique line $O$, called
the \textbf{ortholine} of the trilateral. The ortholine $O$ of the dual
trilateral of a triangle is dual to the orthocenter $o$ of the triangle.%
%TCIMACRO{\FRAME{fhFU}{3.5639in}{1.774in}{0pt}{\Qcb{Orthocenter $o$ and
%ortholine $O$ of a triangle $\overline{a_{1}a_{2}a_{3}}$ and its dual
%trilateral }}{\Qlb{OrthocenterBoth}}{orthocenterboth.eps}%
%width 3.5639in;  height 1.774in;  depth 0pt;  original-width 5.9169in;
%original-height 2.9313in;  cropleft "0";  croptop "1";  cropright "1";
%cropbottom "0";  filename 'OrthocenterBoth.eps';file-properties "XNPEU";}} }%
%BeginExpansion
\begin{figure}[h]%
\centering
\includegraphics[
height=1.774in,
width=3.5639in
]%
{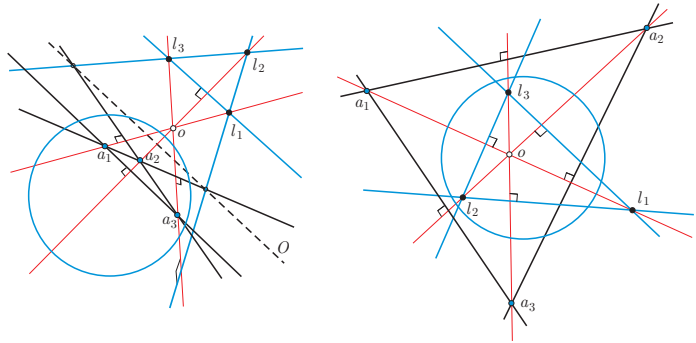}%
\caption{Orthocenter $o$ and ortholine $O$ of a triangle $\overline{a_{1}%
a_{2}a_{3}}$ and its dual trilateral }%
\label{OrthocenterBoth}%
\end{figure}
%EndExpansion

\end{theorem}

Figure \ref{OrthocenterBoth} shows a triangle $\overline{a_{1}a_{2}a_{3}}$,
its dual trilateral with associated triangle $\overline{l_{1}l_{2}l_{3}},$ and
the corresponding orthocenter $o$ and ortholine $O.$

\section{Null points and lines}

In modern treatments of hyperbolic geometry, points at infinity have a
somewhat ambiguous role, and what we call null lines are rarely discussed,
because they are essentially invisible in the Beltrami Poincar\'{e} models
built from differential geometry. However earlier generations of classical
geometers were well aware of them (see for example \cite{Sommerville}).

In universal hyperbolic geometry null points and null lines play a
particularly interesting and important role. Here is a first example, whose
name comes from the definition that a triangle is \textbf{triply nil}
precisely when each of its points is null.

\begin{theorem}
[Triply nil altitudes]Suppose that $\alpha_{1},\alpha_{2}$ and $\alpha_{3}$
are distinct null points, with $b$ any point lying on $\alpha_{1}\alpha_{2}$.
Then the altitude lines to $\alpha_{1}\alpha_{3}$ and $\alpha_{2}\alpha_{3}$
through $b$ are perpendicular.%
%TCIMACRO{\FRAME{fhFU}{3.8622in}{1.6638in}{0pt}{\Qcb{The Triply nil altitudes
%theorem: $bl_{1}\bot bl_{2}$}}{\Qlb{Triply Nil Altitudes}}%
%{niltriangleperpsidealtitudes.eps}{\special{ language "Scientific Word";
%type "GRAPHIC";  maintain-aspect-ratio TRUE;  display "USEDEF";
%valid_file "F";  width 3.8622in;  height 1.6638in;  depth 0pt;
%original-width 6.4539in;  original-height 2.7666in;  cropleft "0";
%filename 'NilTrianglePerpSideAltitudes.eps';file-properties "XNPEU";}} }%
%BeginExpansion
\begin{figure}[h]%
\centering
\includegraphics[
height=1.6638in,
width=3.8622in
]%
{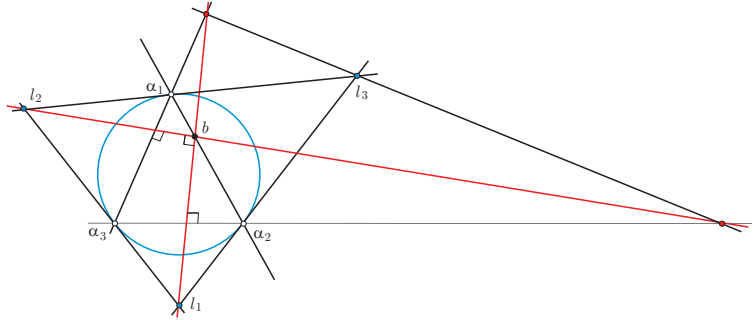}%
\caption{The Triply nil altitudes theorem: $bl_{1}\bot bl_{2}$}%
\label{Triply Nil Altitudes}%
\end{figure}
%EndExpansion

\end{theorem}

Figure \ref{Triply Nil Altitudes} shows that the Triply nil altitudes theorem
may be recast in projective terms: if $l_{1},l_{2}$ and $l_{3}$ are the poles
of the lines of the triangle $\overline{\alpha_{1}\alpha_{2}\alpha_{3}}$ with
respect to the conic $c,$ then the points $\left(  \alpha_{1}\alpha
_{3}\right)  \left(  bl_{1}\right)  =\left(  bl_{2}\right)  ^{\bot},$ $\left(
\alpha_{2}\alpha_{3}\right)  \left(  bl_{2}\right)  =\left(  bl_{1}\right)
^{\bot}$ and $l_{3}$ are collinear.

In keeping with triangles and trilaterals, a (cyclic) set of four points is
called a \textbf{quadrangle}, and a (cyclic) set of four lines a
\textbf{quadrilateral}. The next result restates some facts that we already
know about polarity of cyclic quadrilaterals in terms of perpendicularity.

\begin{theorem}
[Nil quadrangle diagonal]Suppose that $\alpha_{1},\alpha_{2},\alpha_{3}$ and
$\alpha_{4}$ are distinct null points, with diagonal points $e\equiv\left(
\alpha_{1}\alpha_{2}\right)  \left(  \alpha_{3}\alpha_{4}\right)  ,$
$f\equiv\left(  \alpha_{1}\alpha_{3}\right)  \left(  \alpha_{2}\alpha
_{4}\right)  $ and $g\equiv\left(  \alpha_{1}\alpha_{4}\right)  \left(
\alpha_{2}\alpha_{3}\right)  $. Then the lines $ef,$ $eg$ and $fg$ are
mutually perpendicular, and the points $e,f$ and $g$ are also mutually
perpendicular.%
%TCIMACRO{\FRAME{fhFU}{5.6978cm}{4.285cm}{0pt}{\Qcb{Triply right diagonal
%triangle $\overline{efg}$}}{\Qlb{Quad Triply Right}}{nullquadtriplyright.eps}%
%width 5.6978cm;  height 4.285cm;  depth 0pt;  original-width 6.5823cm;
%original-height 4.9291cm;  cropleft "0";  croptop "1";  cropright "1";
%cropbottom "0";  filename 'NullQuadTriplyRight.eps';file-properties "XNPEU";}}
%}%
%BeginExpansion
\begin{figure}[h]%
\centering
\includegraphics[
height=4.285cm,
width=5.6978cm
]%
{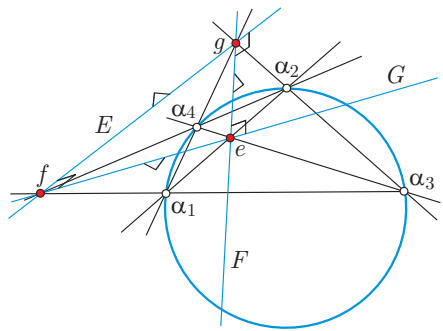}%
\caption{Triply right diagonal triangle $\overline{efg}$}%
\label{Quad Triply Right}%
\end{figure}
%EndExpansion

\end{theorem}

In classical hyperbolic geometry the null points $\alpha_{1},\alpha_{2}%
,\alpha_{3}$ and $\alpha_{4}$ would be considered to be `at infinity', while
the external diagonal points $f$ and $g$ in Figure \ref{Quad Triply Right}
would be invisible. Let us repeat: \textit{for us internal and external points
are equally interesting.}

\section{Couples, parallels and bases}

Points and lines are the basic objects in planar hyperbolic geometry. Given
one point $a,$ we may construct its dual line $A\equiv a^{\perp},$ and
conversely given a line $A$ we may construct its dual point $a\equiv A^{\perp
}.$ After that, there is in general nothing more to construct. For two
objects, namely couples, sides and vertices, the situation is considerably
more interesting, and gives us a chance to introduce some important additional concepts.

Given a non-dual couple $\overline{aL}$ we know we can construct the dual line
$A\equiv a^{\perp}$ and the dual point $l\equiv L^{\perp},$ and the altitude
line $N$ and the altitude point $n.$ Now we introduce another major point of
departure from classical hyperbolic geometry, which provides an ironic twist
to the oft-repeated history of hyperbolic geometry as a development arising
from Euclid's Parallel Postulate.

\begin{theorem}
[Parallel line and point]For any non-dual couple $\overline{aL}$ there is a
unique line $P$ which passes through $a$ and is perpendicular to the altitude
line $N$ of $\overline{aL}$, namely $P\equiv a\left(  a^{\bot}L\right)  $, and
there is a unique point $p$ which lies on $a^{\perp}$ and is perpendicular to
the altitude point $n$ of $\overline{aL}$, namely $p\equiv a^{\bot}\left(
aL^{\bot}\right)  $. Furthermore $P$ and $p$ are dual.%
%TCIMACRO{\FRAME{fhFU}{3.105in}{1.5766in}{0pt}{\Qcb{Parallel line $P$ and
%parallel point $p$ of the couple $\overline{aL}$}}{\Qlb{ParallelLinePoint}%
%}{parallelsboth.eps}{\special{ language "Scientific Word";  type "GRAPHIC";
%width 3.105in;  height 1.5766in;  depth 0pt;  original-width 5.6031in;
%original-height 2.8319in;  cropleft "0";  croptop "1";  cropright "1";
%cropbottom "0";  filename 'ParallelsBoth.EPS';file-properties "XNPEU";}} }%
%BeginExpansion
\begin{figure}[h]%
\centering
\includegraphics[
height=1.5766in,
width=3.105in
]%
{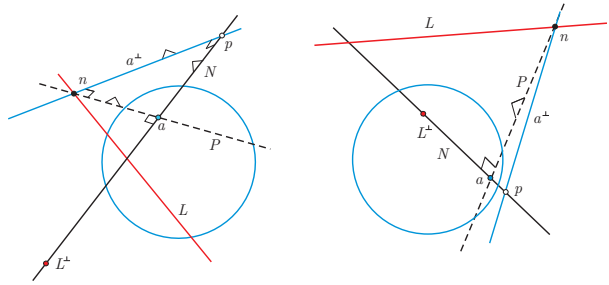}%
\caption{Parallel line $P$ and parallel point $p$ of the couple $\overline
{aL}$}%
\label{ParallelLinePoint}%
\end{figure}
%EndExpansion

\end{theorem}

The line $P$ is the \textbf{parallel line }of the couple $\overline{aL},$ and
the point $p$ is the \textbf{parallel point }of the couple $\overline{aL} $.
These are shown in Figure \ref{ParallelLinePoint}. We may also refer to $P $
as the \textbf{parallel line to the line }$L$\textbf{\ through }$a.$ This is
how we will henceforth use the term \textit{parallel }in universal hyperbolic
geometry---we do not say that \textit{two lines are parallel}.

\begin{theorem}
[Base point and line]For any non-dual couple $\overline{aL}$ there is a unique
point $b$ which lies on both $L$ and the altitude line $N$ of $\overline{aL}$,
namely $b\equiv\left(  aL^{\bot}\right)  L$, and there is a unique line $B$
which passes through both $L^{\bot}$ and the altitude point $n$ of
$\overline{aL}$, namely $B\equiv\left(  a^{\bot}L\right)  L^{\bot}.$
Furthermore $B$ and $b $ are dual.
\end{theorem}

\begin{figure}[h]%
\centering
\includegraphics[
height=1.7036in,
width=3.3549in
]%
{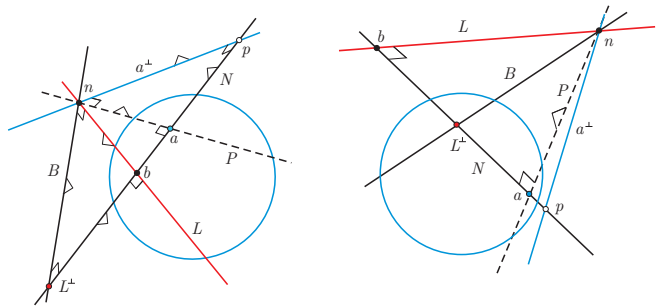}%
\caption{Base point $b$ and base line $B$ of the couple $\overline{aL}$}%
\label{BasePoint}%
\end{figure}

The point $b$ is the \textbf{base point} of the couple $\overline{aL},$ and
the line $B$ is the \textbf{base line }of the couple $\overline{aL}.$ These
are shown in Figure \ref{BasePoint}.%
%TCIMACRO{\FRAME{fhFU}{3.3549in}{1.7036in}{0pt}{\Qcb{Base point $b$ and base
%line $B$ of the couple $\overline{aL}$}}{\Qlb{BasePoint}}{basepointboth.EPS}%
%width 3.3549in;  height 1.7036in;  depth 0pt;  original-width 5.6031in;
%original-height 2.8319in;  cropleft "0";  croptop "1";  cropright "1";
%cropbottom "0";  filename 'BasePointBoth.EPS';file-properties "XNPEU";}} }%
%BeginExpansion
%EndExpansion

\section{Sides, vertices and conjugates}

There are also constructions for sides and vertices, and in fact the complete
pictures associated with the general couple, side and vertex are essentially
the same, but with different objects playing different roles.

\begin{figure}[h]%
\centering
\includegraphics[
height=1.4413in,
width=2.386in
]%
{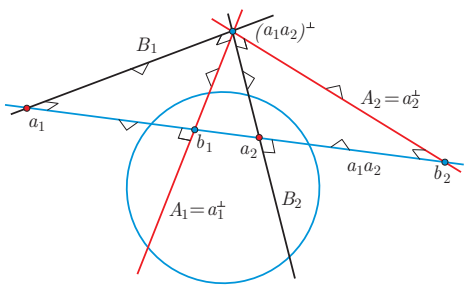}%
\caption{Conjugate points and conjugate lines of the side $\overline
{a_{1}a_{2}}$}%
\label{SideConjugates}%
\end{figure}

Figure \ref{SideConjugates} shows the constructions possible if we start with
a side $\overline{a_{1}a_{2}}.$ Note that in this example one of these points
is internal and one is external. A bit more terminology: we say that the side
$\overline{a_{1}a_{2}}$ is \textbf{null} precisely when $a_{1}a_{2}$ is a null
line, and \textbf{nil }precisely when at least one of $a_{1}$, $a_{2}$ is a
null point. The vertex $\overline{L_{1}L_{2}}$ is \textbf{null} precisely when
$L_{1}L_{2}$ is a null point, and \textbf{nil }precisely when at least one of
$L_{1}$, $L_{2}$ is a null line.

\begin{theorem}
[Side conjugate points and lines ]For any side $\overline{a_{1}a_{2}}$ which
is not both nil and null, there is a unique point $b_{1}\equiv\left(
a_{1}a_{2}\right)  a_{1}^{\bot}$ which lies on $a_{1}a_{2}$ and is
perpendicular to $a_{1},$ and there is a unique point $b_{2}\equiv\left(
a_{1}a_{2}\right)  a_{1}^{\bot}$ which lies on $a_{1}a_{2}$ and is
perpendicular to $a_{2}.$
\end{theorem}

The points $b_{1}$ and $b_{2}$ are the \textbf{conjugate points }of the side
$\overline{a_{1}a_{2}}.$ The duals of these points are the lines $B_{1}\equiv
a_{1}\left(  a_{1}a_{2}\right)  ^{\bot}$ and $B_{2}\equiv a_{2}\left(
a_{1}a_{2}\right)  ^{\bot}$, which are the \textbf{conjugate lines} of the
side $\overline{a_{1}a_{2}}.$ These relations are involutory: if the side
$\overline{b_{1}b_{2}}$ is conjugate to the side $\overline{a_{1}a_{2}},$ then
$\overline{a_{1}a_{2}}$ is also conjugate to $\overline{b_{1}b_{2}}.$
%TCIMACRO{\FRAME{fhFU}{2.386in}{1.4413in}{0pt}{\Qcb{Conjugate points and
%conjugate lines of the side $\overline{a_{1}a_{2}}$}}{\Qlb{SideConjugates}%
%}{conjugatepoints.EPS}{\special{ language "Scientific Word";  type "GRAPHIC";
%width 2.386in;  height 1.4413in;  depth 0pt;  original-width 3.428in;
%original-height 2.0614in;  cropleft "0";  croptop "1";  cropright "1";
%cropbottom "0";  filename 'ConjugatePoints.EPS';file-properties "XNPEU";}} }%
%BeginExpansion
%EndExpansion

The same picture can also be reinterpreted by starting with a vertex.

\begin{theorem}
[Vertex conjugate points and lines]For any vertex $\overline{A_{1}A_{2}}$
which is not both nil and null, there is a unique line $B_{1}\equiv\left(
A_{1}A_{2}\right)  A_{1}^{\bot}$ which passes through $A_{1}A_{2}$ and is
perpendicular to $A_{1},$ and there is a unique line $B_{2}\equiv\left(
A_{1}A_{2}\right)  A_{1}^{\bot}$ which passes through $A_{1}A_{2}$ and is
perpendicular to $A_{2}.$
\end{theorem}

\begin{figure}[h]%
\centering
\includegraphics[
height=1.5351in,
width=2.4608in
]%
{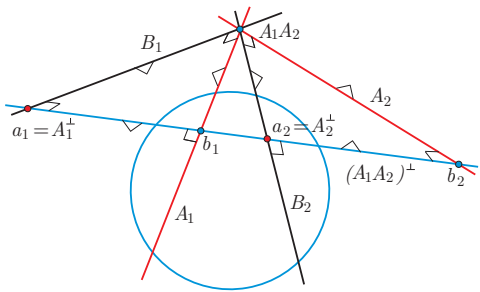}%
\caption{Conjugate points and conjugate lines of the vertex $\overline
{A_{1}A_{2}}$}%
\label{VertexConjugates}%
\end{figure}

The lines $B_{1}$ and $B_{2}$ are the \textbf{conjugate lines }of the vertex
$\overline{A_{1}A_{2}}.$ The duals of these lines are the points $b_{1}\equiv
A_{1}\left(  A_{1}A_{2}\right)  ^{\bot}$ and $b_{2}\equiv A_{2}\left(
A_{1}A_{2}\right)  ^{\bot}$, which are the \textbf{conjugate points }of the
vertex $\overline{A_{1}A_{2}}.$ This relation is also involutory: if the
vertex $\overline{B_{1}B_{2}}$ is conjugate to the vertex $\overline
{A_{1}A_{2}},$ then $\overline{A_{1}A_{2}}$ is also conjugate to
$\overline{B_{1}B_{2}}.$ This is shown also in Figure \ref{VertexConjugates},
which is essentially the same as Figures \ref{BasePoint} and
\ref{SideConjugates}.%
%TCIMACRO{\FRAME{fhFU}{2.4608in}{1.5351in}{0pt}{\Qcb{Conjugate points and
%conjugate lines of the vertex $\overline{A_{1}A_{2}}$}}{\Qlb{VertexConjugates}%
%}{vertexconjugates.EPS}{\special{ language "Scientific Word";
%type "GRAPHIC";  maintain-aspect-ratio TRUE;  display "USEDEF";
%valid_file "F";  width 2.4608in;  height 1.5351in;  depth 0pt;
%original-width 3.428in;  original-height 2.1295in;  cropleft "0";
%filename 'VertexConjugates.EPS';file-properties "XNPEU";}} }%
%BeginExpansion
%EndExpansion

Furthermore the exact same configuration results if we had started with one of
the sides $\overline{a_{1}b_{2}},$ or $\overline{a_{2}b_{1}},$ or
$\overline{b_{1}b_{2}},$ or with one of the vertices $\overline{A_{1}B_{2}},$
or $\overline{A_{2}B_{1}},$ or $\overline{B_{1}B_{2}}.$ However if we had
started with a \textbf{right side}, meaning that the two points are
perpendicular, such as $\overline{a_{1}b_{1}},$ then we would obtain
$L=a_{1}b_{1}$ and its dual point $l.$ But then the conjugate side of
$\overline{a_{1}b_{1}}$ would coincide with $\overline{a_{1}b_{1}},$ and
similarly the conjugate vertex of a right vertex such as $\overline{A_{1}%
B_{1}}$ would coincide with $\overline{A_{1}B_{1}}.$

\section{Reflections}

The basic symmetries of hyperbolic geometry are \textit{reflections, }but they
have a somewhat different character from Euclidean reflections. Hyperbolic
reflections send points to points and lines to lines, preserving incidence, in
other words they are \textit{projective transformations}. There are two
seemingly different notions, the reflection $\sigma_{a}$ in a (non-null) point
$a$, and the reflection $\sigma_{L}$ in a (non-null) line $L$. It is an
important fact that these two notions end up agreeing, in the sense that%
\[
\sigma_{a}=\sigma_{A}%
\]
when $A=a^{\bot}.$

The transformation $\sigma_{a}$ is defined first by its action on null points,
and then by its action on more general points and lines.

For a non-null point $a,$ the reflection $\sigma_{a}$ sends a null point
$\alpha$ to the other null point $\alpha^{\prime}$ on the line $\alpha a.$ We
write
\[
\alpha^{\prime}=\alpha\sigma_{a}.
\]
In case $a\alpha$ is a null line, in other words a tangent to the null circle
$c,$ then $\alpha^{\prime}\equiv\alpha.$ Note that if $a$ was itself a null
point, then this definition would yield a transformation that would send every
null point to $a,$ which will not be a symmetry in the sense we wish. However
such transformations can still be useful.%
%TCIMACRO{\FRAME{fhFU}{3.0328in}{1.6314in}{0pt}{\Qcb{Reflection in the point
%$a$ or the dual line $A\equiv a^{\bot}$}}{\Qlb{Reflection Def}}%
%{reflectionpointboth.eps}{\special{ language "Scientific Word";
%type "GRAPHIC";  maintain-aspect-ratio TRUE;  display "USEDEF";
%valid_file "F";  width 3.0328in;  height 1.6314in;  depth 0pt;
%original-width 4.9475in;  original-height 2.6496in;  cropleft "0";
%filename 'ReflectionPointBoth.eps';file-properties "XNPEU";}} }%
%BeginExpansion
\begin{figure}[h]%
\centering
\includegraphics[
height=1.6314in,
width=3.0328in
]%
{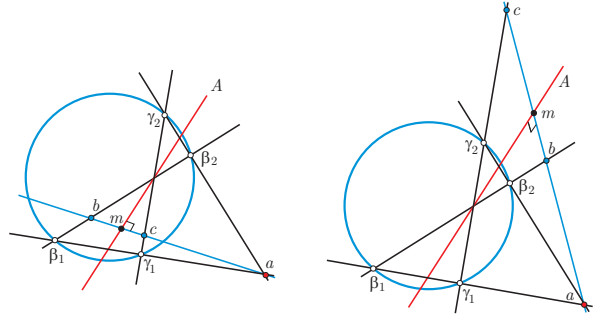}%
\caption{Reflection in the point $a$ or the dual line $A\equiv a^{\bot}$}%
\label{Reflection Def}%
\end{figure}
%EndExpansion

Once the action of a projective transformation on null points is known, it is
determined on all points and lines, first of all on lines through two null
points, and then on an arbitrary point by means of two such lines passing
through it, and then on arbitrary lines. Figure \ref{Reflection Def} shows the
reflection $\sigma_{a}$ and its action on a point $b$ to get $c\equiv
b\sigma_{a}.$ First find a line through $b$ meeting the null circle at points
$\beta_{1}$ and $\beta_{2}$, then construct $\gamma_{1}\equiv\beta_{1}%
\sigma_{a}$ and $\gamma_{1}\equiv\beta_{1}\sigma_{a},$ and then set
\[
c\equiv b\sigma_{a}\equiv\left(  ab\right)  \left(  \gamma_{1}\gamma
_{2}\right)  .
\]

For the reflection $\sigma_{L}$ in a line $L$, the idea is dual to the above.
It is defined first by its action on null lines, and then by its action on
more general points and lines. For a non-null line $L,$ the reflection
$\sigma_{L}$ sends a null line $\Pi$ to the other null line $\Pi^{\prime}$
passing through the point $L\Pi$. We write
\[
\Pi^{\prime}=\Pi\sigma_{L}.
\]
In case $L\Pi$ is a null point, $\Pi^{\prime}\equiv\Pi.$

Once the action of a projective transformation on null lines is known, it is
determined on all points and lines, since it is first of all determined on
points lying on two null lines, and then on an arbitrary line by means of two
such points lying on it, and then on arbitrary points.

Of course there is also a linear algebra/matrix approach to defining
reflections. If $a=\left[  u:v:w\right]  $ then the action of $\sigma
_{a}=\sigma_{a^{\perp}}$ on a point $b\equiv\left[  x:y:z\right]  $ is given
by the projective matrix product%
\[
b\sigma_{a}=\left[  x:y:z\right]
\begin{bmatrix}
u^{2}-v^{2}+w^{2} & 2uv & 2uw\\
2uv & -u^{2}+v^{2}+w^{2} & 2vw\\
-2uw & -2vw & u^{2}+v^{2}-w^{2}%
\end{bmatrix}
\]
where the entries are only determined up to a scalar. The move to three
dimensions simplifies the discussion.

\section{Midpoints, midlines, bilines and bipoints}

The notion of the \textit{midpoint }of a side can be defined once we have the
notion of a reflection. It also has a metrical formulation in terms of
quadrance, which we have not yet introduced. There are three other closely
related concepts, that of \textit{midline}, \textit{biline} and
\textit{bipoint}. Midpoints and midlines refer to sides, while bilines and
bipoints refer to vertices. The existence of these objects reduces to
questions in number theory----whether or not certain quadratic equations have solutions.

If $c=b\sigma_{d}$ then we say $d$ is a\textbf{\ midpoint} of the side
$\overline{bc}.$ In this case the point $e\equiv d^{\bot}\left(  bc\right)  $
is also a midpoint of $\overline{bc}$, and the two midpoints $d$ and $e$ of
$\overline{bc}$ are perpendicular. The dual line $D\equiv d^{\bot}$ is a
\textbf{midline }of $\overline{bc},$ meaning that it meets $bc$
perpendicularly in a midpoint, namely $e$. The reflection $\sigma_{D}$ in $D $
of $b$ is also $c.$ Similarly the dual line $E\equiv e^{\bot}$ is also a
midline of $\overline{bc}.$ In Euclidean geometry midlines are called
\textit{perpendicular bisectors}. In hyperbolic geometry there are generally
either zero or two midpoints between any two points, and so also zero or two midlines.

\begin{figure}[h]%
\centering
\includegraphics[
height=1.7991in,
width=3.6969in
]%
{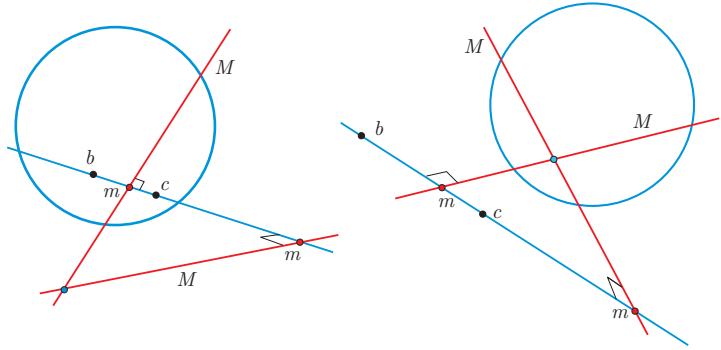}%
\caption{Midpoints and midlines of $\overline{bc}$}%
\label{MidpointsMidlines}%
\end{figure}

\begin{figure}[h]%
\centering
\includegraphics[
height=1.8755in,
width=3.1283in
]%
{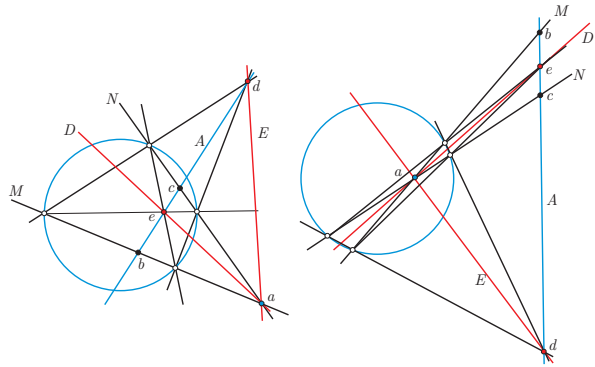}%
\caption{Midpoints $d$ and $e$ of $\overline{bc},$ or bilines $D$ and $E$ of
$\overline{MN}$ }%
\label{Midpoints Bisectors}%
\end{figure}

Figure \ref{MidpointsMidlines} shows the midpoints $m$ and midlines $M$ of a
side $\overline{bc}.$
%TCIMACRO{\FRAME{fhFU}{3.6969in}{1.7991in}{0pt}{\Qcb{Midpoints and midlines of
%$\overline{bc}$}}{\Qlb{MidpointsMidlines}}{midpoints.eps}%
%width 3.6969in;  height 1.7991in;  depth 0pt;  original-width 5.348in;
%original-height 2.5884in;  cropleft "0";  croptop "1";  cropright "1";
%cropbottom "0";  filename 'Midpoints.eps';file-properties "XNPEU";}} }%
%BeginExpansion
%EndExpansion
Figure \ref{Midpoints Bisectors} shows how to construct the midpoints and
midlines of the side $\overline{bc}$ when such exist. For the midpoints of
$\overline{bc}$ we first join $b$ and $c$ to the point $a\equiv\left(
bc\right)  ^{\bot}$ to form lines $M$ and $N.$ If these are both interior
lines, then their meets with the null circle give a completely nil quadrangle
one of whose diagonal points is $a,$ and the other two diagonal points $d$ and
$e$ lie on $bc$ and are the required midpoints. The duals $D$ and $E$ of $d$
and $e$ respectively are the midlines of $\overline{bc}.$
%TCIMACRO{\FRAME{fhFU}{3.1283in}{1.8755in}{0pt}{\Qcb{Midpoints $d$ and $e$ of
%$\overline{bc},$ or bilines $D$ and $E$ of $\overline{MN}$ }}%
%{\Qlb{Midpoints Bisectors}}{BisectorsBoth.eps}%
%width 3.1283in;  height 1.8755in;  depth 0pt;  original-width 5.6361in;
%original-height 3.3678in;  cropleft "0";  croptop "1";  cropright "1";
%cropbottom "0";  filename 'BisectorsBoth.eps';file-properties "XNPEU";}} }%
%BeginExpansion
%EndExpansion

The dual notions to midpoints and midlines of sides are the notions
of\ bilines and bipoints of vertices. If $M$ and $N$ are lines and
$M=N\sigma_{D}$ then the line $D$ is a \textbf{biline }of the vertex
$\overline{MN}.$ In this case the line $E\equiv D^{\bot}\left(  MN\right)  $
is also a biline of $\overline{MN},$ and the two bilines $D$ and $E$ of
$\overline{MN}$ are perpendicular.

The dual point $d\equiv D^{\bot}$ is a \textbf{bipoint} of $\overline{MN},$
and it joins $MN$ perpendicularly in a biline, namely $E.$ This implies that
the reflection $\sigma_{d}$ in $d$ of $M$ is also $N.$ Similarly the dual
point $e\equiv E^{\bot}$ is also a bipoint of $\overline{MN}.$ In Euclidean
geometry, bilines are called vertex bisectors or angle bisectors. Bipoints
have no Euclidean analog.

Figure \ref{Midpoints Bisectors} can equally well be interpreted as
illustrating the process of obtaining bilines $D$ and $E$ and bipoints $d$ and
$e$ of the vertex $\overline{MN}.$

\section{Hyperbolic triangle geometry}

The richness of Euclidean triangle geometry is not reflected in the classical
hyperbolic setting, but the situation is remedied with universal hyperbolic
geometry. Here we give just a quick glimpse in this fascinating direction,
which will be the focus of a subsequent paper in this series.%
%TCIMACRO{\FRAME{fhFU}{3.3076in}{2.1775in}{0pt}{\Qcb{Circumcenters and
%circumlines of a triangle $\overline{a_{1}a_{2}a_{3}}$}}%
%{\Qlb{CircumcentersCircumlines}}{circumcenterscircumlines.eps}%
%width 3.3076in;  height 2.1775in;  depth 0pt;  original-width 6.1594in;
%original-height 4.044in;  cropleft "0";  croptop "1";  cropright "1";
%filename 'CircumcentersCircumlines.EPS';file-properties "XNPEU";}} }%
%BeginExpansion
\begin{figure}[h]%
\centering
\includegraphics[
height=2.1775in,
width=3.3076in
]%
{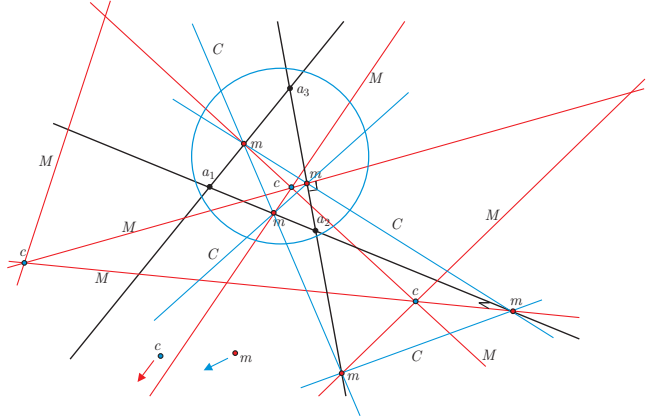}%
\caption{Circumcenters and circumlines of a triangle $\overline{a_{1}%
a_{2}a_{3}}$}%
\label{CircumcentersCircumlines}%
\end{figure}
%EndExpansion

Figure \ref{CircumcentersCircumlines} shows a triangle $\overline{a_{1}%
a_{2}a_{3}}$ together with its six\textit{\ }midpoints $m$ (one is off the
page) and corresponding six\textit{\ }midlines $M$. The midpoints are
collinear three at a time on four lines $C$ called \textbf{circumlines}. The
midlines are concurrent three at a time on four points $c$ called
\textbf{circumcenters}. The circumcenters are dual to the circumlines.
Although we have not defined circles yet, a triangle generally has zero or
four circumcircles, whose centers are at its circumcenters, if these exist.

\begin{figure}[h]%
\centering
\includegraphics[
height=6.7656cm,
width=10.7961cm
]%
{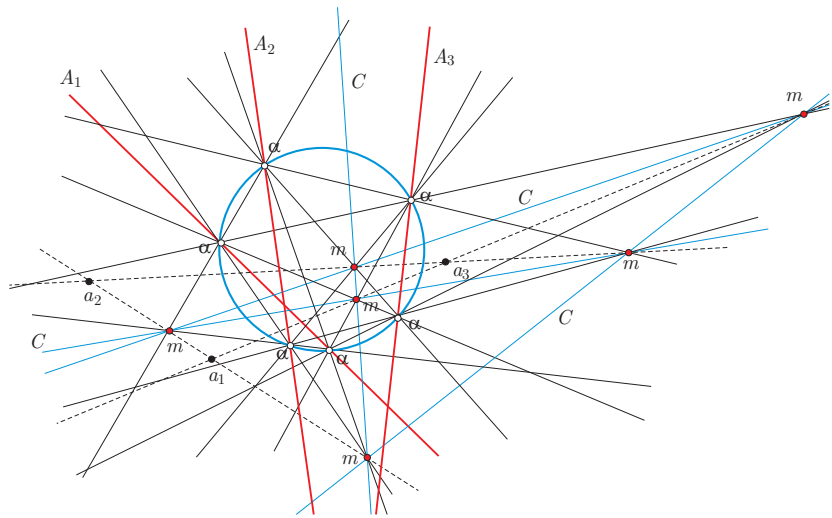}%
\caption{Pascal's theorem via hyperbolic geometry}%
\label{Pascal's theorem}%
\end{figure}

An important application of circumlines\textit{\ }is to \textit{Pascal's
theorem}, one of the great classical results of geometry. Figure
\ref{Pascal's theorem} shows three lines $A_{1},A_{2}$ and $A_{3}$ which meet
the null circle at six null points $\alpha$. The dual points of $A_{1},A_{2} $
and $A_{3}$ are $a_{1},a_{2}$ and $a_{3}$ respectively. The triangle
$\overline{a_{1}a_{2}a_{3}}$ has six midpoints $m$ (red) and the four
circumlines $C$ (blue) pass through three midpoints each.%
%TCIMACRO{\FRAME{fhFU}{10.7961cm}{6.7656cm}{0pt}{\Qcb{Pascal's theorem via
%hyperbolic geometry}}{\Qlb{Pascal's theorem}}{pascalhyperbolic.eps}%
%width 10.7961cm;  height 6.7656cm;  depth 0pt;  original-width 15.371cm;
%original-height 9.6023cm;  cropleft "0";  croptop "1";  cropright "1";
%cropbottom "0";  filename 'PascalHyperbolic.eps';file-properties "XNPEU";}} }%
%BeginExpansion
%EndExpansion

Two of the original three lines, such as $A_{1}$ and $A_{2},$ determine four
null points $\alpha$, and the other two diagonal points formed by such a
quadrangle of null points, not including $A_{1}A_{2},$ give two midpoints $m,
$ in this case of the side $\overline{a_{1}a_{2}}.$ So Pascal's theorem is
here seen as a consequence of the fact that the\textit{\ six midpoints of a
triangle are collinear three at a time} forming the circumlines.

The six null points $\alpha$ can be partitioned into three sets of two in $15$
ways. By different choices of the lines $A_{1},A_{2}$ and $A_{3},$ there are
altogether $15$ such diagrams associated to the same six null points, and $60$
possible circumlines $C$ playing the role of Pascal's line. Such a large
configuration has many remarkable features, some of them projective, some of
them metrical.

\section{Parallels and the doubled triangle}

\begin{figure}[h]%
\centering
\includegraphics[
height=1.6462in,
width=2.2591in
]%
{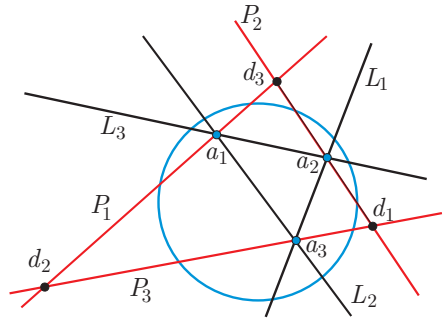}%
\caption{A triangle $\overline{a_{1}a_{2}a_{3}}$ and its double triangle
$\overline{d_{1}d_{2}d_{3}}.$}%
\label{DoubleTriangle}%
\end{figure}

Given a triangle $\overline{a_{1}a_{2}a_{3}},$ the \textbf{double triangle}
$\overline{d_{1}d_{2}d_{3}}$ is the triangle whose lines are the parallels
$P_{1},P_{2}$ and $P_{3}$ to the lines $L_{1},L_{2}$ and $L_{3}$ of
$\overline{a_{1}a_{2}a_{3}}$ through the points $a_{1},a_{2}$ and $a_{3}$
respectively. We retain the usual notational conventions, so that $d_{1}%
=P_{2}P_{3}$ etc. The situation is shown in Figure \ref{DoubleTriangle}.%
%TCIMACRO{\FRAME{fhFU}{2.2591in}{1.6462in}{0pt}{\Qcb{A triangle $\overline
%{a_{1}a_{2}a_{3}}$ and its double triangle $\overline{d_{1}d_{2}d_{3}}.$}%
%}{\Qlb{DoubleTriangle}}{doubletriangle.EPS}%
%width 2.2591in;  height 1.6462in;  depth 0pt;  original-width 2.7789in;
%original-height 2.017in;  cropleft "0";  croptop "1";  cropright "1";
%cropbottom "0";  filename 'DoubleTriangle.EPS';file-properties "XNPEU";}} }%
%BeginExpansion
%EndExpansion

The next theorem is surprising to me, and seems to require a somewhat involved computation.

\begin{theorem}
[Double median triangle]If $\overline{d_{1}d_{2}d_{3}}$ is the double triangle
of a triangle $\overline{a_{1}a_{2}a_{3}},$ then $a_{1},a_{2}$ and $a_{3}$ are
midpoints of the sides of $\overline{d_{1}d_{2}d_{3}}$.
\end{theorem}

In Euclidean geometry the points in the next two theorems would both be the
centroid of the triangle.

\begin{theorem}
[Double point]If $\overline{d_{1}d_{2}d_{3}}$ is the double triangle of a
triangle $\overline{a_{1}a_{2}a_{3}},$ then the lines $a_{1}d_{1}$,
$a_{2}d_{2}$ and $a_{3}d_{3}$ are concurrent in a point $x.$
\end{theorem}

\begin{theorem}
[Second double point]If $\overline{d_{1}d_{2}d_{3}}$ is the double triangle of
a triangle $\overline{a_{1}a_{2}a_{3}},$ and $\overline{g_{1}g_{2}g_{3}}$ is
the double triangle of $\overline{d_{1}d_{2}d_{3}},$ then the lines
$a_{1}g_{1}$, $a_{2}g_{2}$ and $a_{3}g_{3}$ are concurrent in a point $y.$
\end{theorem}

\begin{figure}[h]%
\centering
\includegraphics[
height=1.7684in,
width=2.8941in
]%
{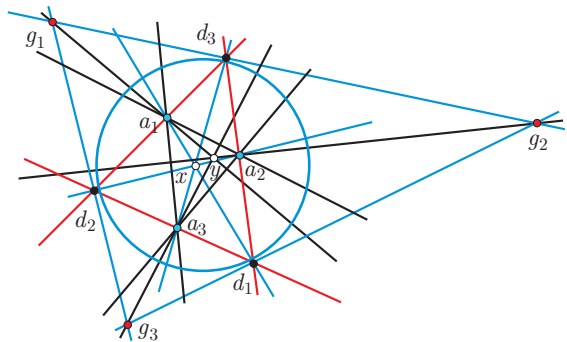}%
\caption{First and second double points of a triangle $\overline{a_{1}%
a_{2}a_{3}}$}%
\label{Double Points}%
\end{figure}

These are shown in Figure \ref{Double Points}; $x$ is the \textbf{double
point} of the triangle $\overline{a_{1}a_{2}a_{3}},$ and $y$ is the
\textbf{second double point} of the triangle $\overline{a_{1}a_{2}a_{3}}.$
%TCIMACRO{\FRAME{fhFU}{2.8941in}{1.7684in}{0pt}{\Qcb{First and second double
%points of a triangle $\overline{a_{1}a_{2}a_{3}}$}}{\Qlb{Double Points}%
%}{doublepoints.EPS}{\special{ language "Scientific Word";  type "GRAPHIC";
%width 2.8941in;  height 1.7684in;  depth 0pt;  original-width 3.6473in;
%original-height 2.2183in;  cropleft "0";  croptop "1";  cropright "1";
%cropbottom "0";  filename 'DoublePoints.EPS';file-properties "XNPEU";}} }%
%BeginExpansion
%EndExpansion

It is not the case that the pattern continues in the obvious way: one cannot
define a third double point in an analogous way. The study of the double
triangle is clearly an interesting departure point from Euclidean triangle geometry.

\section{Quadrance and spread}

We now introduce the two basic measurements in universal hyperbolic geometry,
the \textit{quadrance }between points and the \textit{spread} between lines.
These are analogs of the corresponding notions in Euclidean \textit{rational
trigonometry}, but we assume no familiarity with this theory (although for a
deeper understanding one should carefully compare the two). Our definition of
quadrance and spread follows our projective orientation, and is given here in
terms of the \textit{cross-ratio} between four particular points or lines. The
importance of this particular cross-ratio was already recognized in
(\cite{Brauner}).

\begin{figure}[h]%
\centering
\includegraphics[
height=4.131cm,
width=6.7754cm
]%
{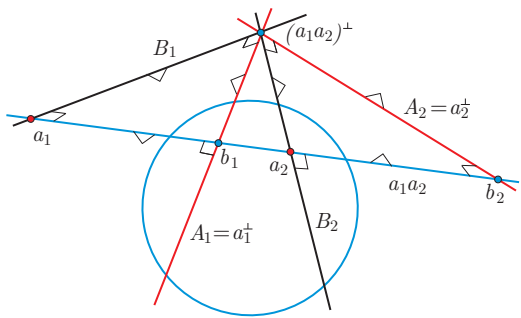}%
\caption{Quadrance defined as cross-ratio: $q\left(  a_{1},a_{2}\right)
\equiv\left(  a_{1},b_{2}:a_{2},b_{1}\right)  $}%
\label{Conjugate points}%
\end{figure}

Suppose that $a_{1}$ and $a_{2}$ are points and that $b_{1}$ and $b_{2}$ are
the conjugate points of the side $\overline{a_{1}a_{2}},$ as shown in Figure
\ref{Conjugate points}. Then define the \textbf{quadrance} between $a_{1}$ and
$a_{2}$ to be the cross-ratio of points:%
\[
q\left(  a_{1},a_{2}\right)  \equiv\left(  a_{1},b_{2}:a_{2},b_{1}\right)  .
\]
The quadrance $q\left(  a_{1},a_{2}\right)  $ is zero if $a_{1}=a_{2}.$ It is
negative if $a_{1}$ and $a_{2}$ are both interior points, and approaches
infinity as $a_{1}$ or $a_{2}$ approaches the null circle. It is undefined (or
infinite) if one or both of $a_{1},a_{2}$ is a null point. It is positive if
one of $a_{1}$ and $a_{2}$ is an interior point and the other is an exterior
point. It is negative if both $a_{1}$ and $a_{2}$ are exterior points and
$a_{1}a_{2}$ is an interior line. It is zero if $a_{1}a_{2}$ is a null line.
It is positive if both $a_{1}$ and $a_{2}$ are exterior points and $a_{1}%
a_{2}$ is an exterior line.%
%TCIMACRO{\FRAME{fhFU}{6.7754cm}{4.131cm}{0pt}{\Qcb{Quadrance defined as
%cross-ratio: $q\left(  a_{1},a_{2}\right)  \equiv\left(  a_{1},b_{2}%
%:a_{2},b_{1}\right)  $}}{\Qlb{Conjugate points}}{sideends.eps}%
%width 6.7754cm;  height 4.131cm;  depth 0pt;  original-width 9.305cm;
%original-height 6.2634cm;  cropleft "0";  croptop "1";  cropright "1";
%cropbottom "0";  filename 'SideEnds.eps';file-properties "XNPEU";}} }%
%BeginExpansion
%EndExpansion

The usual distance $d\left(  a_{1},a_{2}\right)  $ between $a_{1}$ and $a_{2}$
in the Klein model is also defined in terms of a cross-ratio, involving two
\textit{other} points: the meets of $a_{1}a_{2}$ with the null circle. This is
problematic in three ways. First of all for two general points there may be no
such meet, and so the Klein distance does not extend to general points. But
even if the meets exist, it is not easy to separate them algebraically to get
four points in a prescribed and canonical order to apply the cross-ratio.
Finally to get a quantity that acts somewhat linearly, one is forced to
introduce a logarithm or inverse circular function. This is much more
complicated analytically, and makes extending the theory to finite fields, for
example, more problematic.

In any case it turns out that if $a_{1}$ and $a_{2}$ are interior points,
there is a relation between quadrance and the Klein distance:
\begin{equation}
q\left(  a_{1},a_{2}\right)  =-\sinh^{2}\left(  d\left(  a_{1},a_{2}\right)
\right)  . \label{QuadranceDistance}%
\end{equation}

To define the spread between lines, we proceed in a dual fashion. In Figure
\ref{Conjugate points}, $B_{1}$ and $B_{2}$ are the conjugate lines of the
vertex $\overline{A_{1}A_{2}}.$ Define the \textbf{spread }between the lines
$A_{1}$ and $A_{2}$ to be the cross-ratio of lines:%
\[
S\left(  A_{1},A_{2}\right)  \equiv\left(  A_{1},B_{2}:A_{2},B_{1}\right)  .
\]
This is positive if $A_{1}$ and $A_{2}$ are both interior lines that meet in
an interior point. In fact if $\theta\left(  A_{1},A_{2}\right)  $ is the
usual angle between $A_{1}$ and $A_{2}$ in the Klein model, then it turns out
that
\begin{equation}
S\left(  A_{1},A_{2}\right)  =\sin^{2}\left(  \theta\left(  A_{1}%
,A_{2}\right)  \right)  . \label{SpreadAngle2}%
\end{equation}

The relations (\ref{QuadranceDistance}) and (\ref{SpreadAngle2}) allow you to
translate the subsequent theorems in this paper to formulas of classical
hyperbolic trigonometry in the special case of interior points and lines.

The spread $S\left(  A_{1},A_{2}\right)  $ between lines $A_{1}$ and $A_{2}$
is equal to the quadrance between the dual points, that is%
\[
S\left(  A_{1},A_{2}\right)  =q\left(  A_{1}^{\bot},A_{2}^{\bot}\right)  .
\]
So the basic duality between points and lines extends to the two fundamental measurements.

A \textbf{circle} is given by an equation of the form $q\left(  x,a\right)
=k$ for some fixed point $a$ called the \textbf{center}, and a number $k$
called the \textbf{quadrance}.

\begin{figure}[h]%
\centering
\includegraphics[
height=1.6579in,
width=3.8373in
]%
{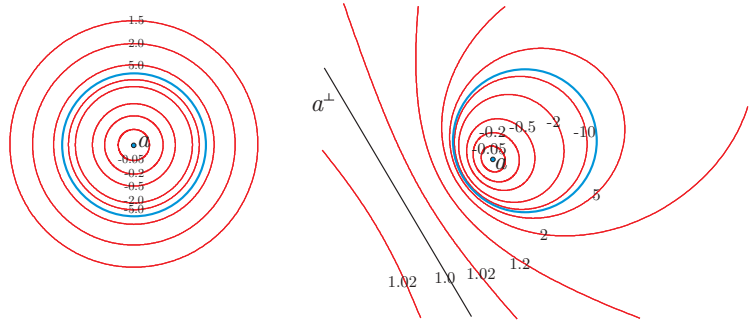}%
\caption{Circles centered at $a$ (interior)}%
\label{CirclesInterior}%
\end{figure}

\begin{figure}[h]%
\centering
\includegraphics[
height=2.0025in,
width=2.704in
]%
{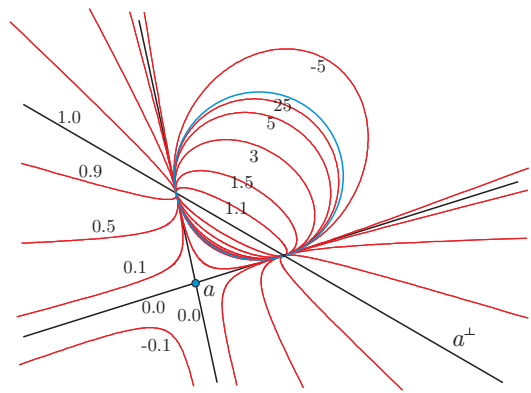}%
\caption{Circles centered at $a$ (exterior)}%
\label{CirclesExterior}%
\end{figure}

Here we show the circles centered at a point $a$ of various quadrances. Figure
\ref{CirclesInterior} shows circles centered at a point $a$ when $a$ is an
interior point. Figure \ref{CirclesExterior} shows circles centered at an
exterior point $a.$ Both of these diagrams should be studied carefully. Note
that the dual line $a^{\perp}$ of $a$ is such a circle, of quadrance $1.$ Also
note that the situation is dramatically different for $a$ interior or $a$
exterior.
%TCIMACRO{\FRAME{fhFU}{3.8373in}{1.6579in}{0pt}{\Qcb{Circles centered at $a$
%(interior)}}{\Qlb{CirclesInterior}}{circlescenterboth.EPS}%
%width 3.8373in;  height 1.6579in;  depth 0pt;  original-width 6.7222in;
%original-height 2.8893in;  cropleft "0";  croptop "1";  cropright "1";
%cropbottom "0";  filename 'CirclesCenterBoth.EPS';file-properties "XNPEU";}}
%}%
%BeginExpansion
%EndExpansion
In the case of $a$ an exterior point, there is a non-trivial circle of
quadrance $0,$ namely the two null lines through $a,$ and all circles meet the
null circle at the two points where the dual line $a^{\bot}$ meets it. In
classical hyperbolic geometry such curves are known as \textit{constant width
curves}---however from our point of view they are just circles.%
%TCIMACRO{\FRAME{fhFU}{2.704in}{2.0025in}{0pt}{\Qcb{Circles centered at $a$
%(exterior)}}{\Qlb{CirclesExterior}}{circlescenterexterior.EPS}%
%width 2.704in;  height 2.0025in;  depth 0pt;  original-width 4.7063in;
%original-height 3.4783in;  cropleft "0";  croptop "1";  cropright "1";
%filename 'CirclesCenterExterior.EPS';file-properties "XNPEU";}} }%
%BeginExpansion
%EndExpansion

There is a dual approach to circles where we use the relation $S\left(
X,L\right)  =k$ for a fixed line $L$ and a variable line $X.$ We leave it to
the reader to show that we obtain the envelope of a circle as defined in terms
of points. So the notion of a circle is essentially self-dual.

\section{Basic trigonometric laws}

For most calculations, we need explicit analytic formulae for the main measurements.

\begin{theorem}
The \textbf{quadrance }between points $a_{1}\equiv\left[  x_{1}:y_{1}%
:z_{1}\right]  $ and $a_{2}\equiv\left[  x_{2}:y_{2}:z_{2}\right]  $ is%
\[
q\left(  a_{1},a_{2}\right)  =1-\frac{\left(  x_{1}x_{2}+y_{1}y_{2}-z_{1}%
z_{2}\right)  ^{2}}{\left(  x_{1}^{2}+y_{1}^{2}-z_{1}^{2}\right)  \left(
x_{2}^{2}+y_{2}^{2}-z_{2}^{2}\right)  }.
\]

\end{theorem}

\begin{theorem}
The \textbf{spread} between lines $L_{1}\equiv\left(  l_{1}:m_{1}%
:n_{1}\right)  $ and $L_{2}\equiv\left(  l_{2}:m_{2}:n_{2}\right)  $ is%
\[
S\left(  L_{1},L_{2}\right)  =1-\frac{\left(  l_{1}l_{2}+m_{1}m_{2}-n_{1}%
n_{2}\right)  ^{2}}{\left(  l_{1}^{2}+m_{1}^{2}-n_{1}^{2}\right)
\allowbreak\left(  l_{2}^{2}+m_{2}^{2}-n_{2}^{2}\right)  }.
\]

\end{theorem}

These expressions are not defined if one or more of the points or lines
involved is null, and reinforce the fact that the duality between points and
lines extends to quadrances and spreads, and so every metrical result can be
expected to have a dual formulation.

For a triangle $\overline{a_{1}a_{2}a_{3}}$ with associated trilateral
$\overline{L_{1}L_{2}L_{3}}$ we will use the usual convention that
$q_{1}\equiv q\left(  a_{2},a_{3}\right)  $, $q_{2}\equiv q\left(  a_{1}%
,a_{3}\right)  $ and $q_{3}\equiv q\left(  a_{1},a_{2}\right)  $, and
$S_{1}\equiv S\left(  L_{2},L_{3}\right)  $, $S_{2}\equiv S\left(  L_{1}%
,L_{3}\right)  $ and $S_{3}\equiv S\left(  L_{1},L_{2}\right)  $. This
notation will also be used in the degenerate case when $a_{1},a_{2}$ and
$a_{3}$ are collinear, or $L_{1},L_{2}$ and $L_{3}$ are concurrent.
%TCIMACRO{\FRAME{fhFU}{5.7063cm}{3.7641cm}{0pt}{\Qcb{Quadrance and spreads in a
%hyperbolic triangle}}{\Qlb{TriNotation0}}{trianglenotation.eps}%
%width 5.7063cm;  height 3.7641cm;  depth 0pt;  original-width 5.8165cm;
%original-height 4.0244cm;  cropleft "0";  croptop "1";  cropright "1";
%cropbottom "0";  filename 'TriangleNotation.eps';file-properties "XNPEU";}} }%
%BeginExpansion
\begin{figure}[h]%
\centering
\includegraphics[
height=3.7641cm,
width=5.7063cm
]%
{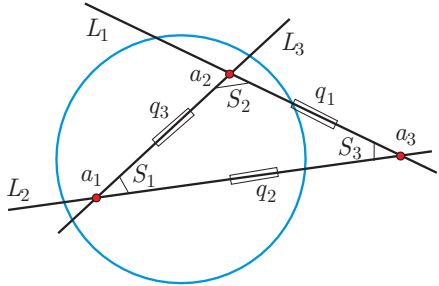}%
\caption{Quadrance and spreads in a hyperbolic triangle}%
\label{TriNotation0}%
\end{figure}
%EndExpansion

With this notation, here are the main trigonometric laws in the subject. These
are among the most important formulas in mathematics.

\begin{theorem}
[Triple quad formula]\textit{If }$a_{1},a_{2}$\textit{\ and }$a_{3}%
$\textit{\ are collinear points then }%
\[
\left(  q_{1}+q_{2}+q_{3}\right)  ^{2}=2\left(  q_{1}^{2}+q_{2}^{2}+q_{3}%
^{2}\right)  +4q_{1}q_{2}q_{3}.
\]

\end{theorem}

\begin{theorem}
[Triple spread formula]\textit{If }$L_{1},L_{2}$\textit{\ and }$L_{3}%
$\textit{\ are concurrent lines then }%
\[
\left(  S_{1}+S_{2}+S_{3}\right)  ^{2}=2\left(  S_{1}^{2}+S_{2}^{2}+S_{3}%
^{2}\right)  +4S_{1}S_{2}S_{3}.
\]

\end{theorem}

\begin{theorem}
[Pythagoras]\textit{If }$L_{1}$\textit{\ and }$L_{2}$\textit{\ are
perpendicular lines then}%
\[
q_{3}=q_{1}+q_{2}-q_{1}q_{2}.
\]

\end{theorem}

\begin{theorem}
[Pythagoras' dual]\textit{If }$a_{1}$ and $a_{2}$ a\textit{re perpendicular
points then}%
\[
S_{3}=S_{1}+S_{2}-S_{1}S_{2}.
\]

\end{theorem}

\begin{theorem}
[Spread law]%
\[
\frac{S_{1}}{q_{1}}=\frac{S_{2}}{q_{2}}=\frac{S_{3}}{q_{3}}.
\]

\end{theorem}

\begin{theorem}
[Spread dual law]%
\[
\frac{q_{1}}{S_{1}}=\frac{q_{2}}{S_{2}}=\frac{q_{3}}{S_{3}}.
\]

\end{theorem}

\begin{theorem}
[Cross law]%
\[
\left(  q_{1}q_{2}S_{3}-\left(  q_{1}+q_{2}+q_{3}\right)  +2\right)
^{2}=4\left(  1-q_{1}\right)  \left(  1-q_{2}\right)  \left(  1-q_{3}\right)
.
\]

\end{theorem}

\begin{theorem}
[Cross dual law]%
\[
\left(  S_{1}S_{2}q_{3}-\left(  S_{1}+S_{2}+S_{3}\right)  +2\right)
^{2}=4\left(  1-S_{1}\right)  \left(  1-S_{2}\right)  \left(  1-S_{3}\right)
.
\]

\end{theorem}

There are three symmetrical forms of Pythagoras' theorem, the Cross law and
their duals, obtained by rotating indices. A proper appreciation for the
beauty and power of these formulas requires some familiarity with rational
trigonometry in the plane (see \cite{Wild1}), together with rolling up one's
sleeves and solving many trigonometric problems in the hyperbolic setting. For
students of geometry, this is an excellent undertaking.

\section{Right triangles and trilaterals}

Right triangles and trilaterals have some additional important properties
besides the fundamental Pythagoras theorem we have already mentioned. We leave
the dual results to the reader. Thales' theorem shows that there is an aspect
of similar triangles in hyperbolic geometry. It also helps explain why spread
is the primary measurement between lines in rational trigonometry.

\begin{theorem}
[Thales]\textit{Suppose that }$\overline{a_{1}a_{2}a_{3}}$\textit{\ is a right
triangle with }$S_{3}=1.$\textit{\ Then }%
\[
S_{1}=\frac{q_{1}}{q_{3}}\qquad\mathrm{and}\qquad S_{2}=\frac{q_{2}}{q_{3}}.
\]

\end{theorem}

%

%TCIMACRO{\FRAME{fhFU}{3.571in}{1.617in}{0pt}{\Qcb{Thales' theorem:
%$S_{1}=q_{1}/q_{3}$}}{\Qlb{Thales}}{thalesboth3.EPS}%
%width 3.571in;  height 1.617in;  depth 0pt;  original-width 4.827in;
%original-height 2.1696in;  cropleft "0";  croptop "1";  cropright "1";
%cropbottom "0";  filename 'ThalesBoth3.EPS';file-properties "XNPEU";}} }%
%BeginExpansion
\begin{figure}[h]%
\centering
\includegraphics[
height=1.617in,
width=3.571in
]%
{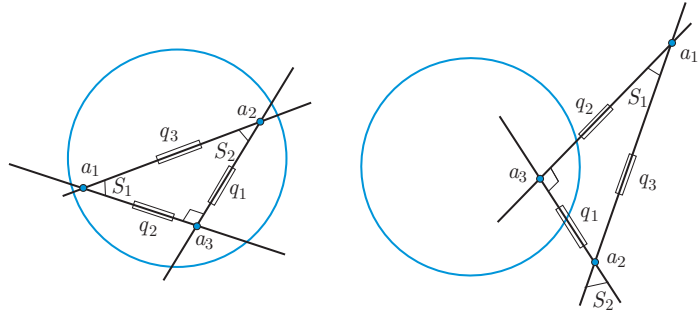}%
\caption{Thales' theorem: $S_{1}=q_{1}/q_{3}$}%
\label{Thales}%
\end{figure}
%EndExpansion

The Right parallax theorem generalizes, and dramatically simplifies, a famous
formula of Bolyai and Lobachevsky (see \cite{Greenberg}) which usually
requires exponential and circular functions, hence a prior understanding of
real numbers.

\begin{theorem}
[Right parallax]\textit{If a right triangle }$\overline{a_{1}a_{2}a_{3}}%
$\textit{\ has spreads }$S_{1}=0$\textit{, }$S_{2}=S$\textit{\ and }$S_{3}%
=1,$\textit{\ then it will have only one defined quadrance }$q_{1}%
=q$\textit{\ given by}%
\[
q=\frac{S-1}{S}.
\]

\end{theorem}

%

%TCIMACRO{\FRAME{fhFU}{3.9502in}{1.4255in}{0pt}{\Qcb{Right parallax
%theorem:\ $q=\left(  S-1\right)  /S$}}{\Qlb{RightParallax}}%
%{rightparallax3.EPS}{\special{ language "Scientific Word";  type "GRAPHIC";
%width 3.9502in;  height 1.4255in;  depth 0pt;  original-width 4.5278in;
%original-height 1.6171in;  cropleft "0";  croptop "1";  cropright "1";
%cropbottom "0";  filename 'RightParallax3.EPS';file-properties "XNPEU";}} }%
%BeginExpansion
\begin{figure}[h]%
\centering
\includegraphics[
height=1.4255in,
width=3.9502in
]%
{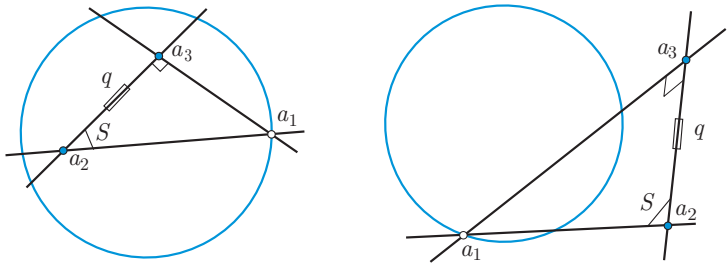}%
\caption{Right parallax theorem:\ $q=\left(  S-1\right)  /S$}%
\label{RightParallax}%
\end{figure}
%EndExpansion
We may restate this result in the form%
\[
S=\frac{1}{1-q}.
\]

Napier's Rules are much simpler in the universal setting, where only high
school algebra is required.

\begin{theorem}
[Napier's Rules]\textit{Suppose a right triangle }$\overline{a_{1}a_{2}a_{3}}%
$\textit{\ has quadrances }$q_{1},q_{2}$\textit{\ and }$q_{3},$\textit{\ and
spreads }$S_{1},S_{2}$\textit{\ and }$S_{3}=1.$\textit{\ Then any two of the
five quantities }$S_{1},S_{2},q_{1},q_{2}$ and $q_{3}$\textit{\ determine the
other three, solely by the three basic equations from Thales' theorem and
Pythagoras' theorem: }%
\[
S_{1}=\frac{q_{1}}{q_{3}}\qquad S_{2}=\frac{q_{2}}{q_{3}}\qquad q_{3}%
=q_{1}+q_{2}-q_{1}q_{2}.
\]

\end{theorem}

\section{Triangle proportions and barycentric coordinates}

The following theorems implicitly involve barycentric coordinates. These are
quite useful both in universal and classical hyperbolic geometry, see for
example \cite{Ungar}.

\begin{theorem}
[Triangle proportions]\textit{Suppose that }$\overline{a_{1}a_{2}a_{3}}%
$\textit{\ is a triangle with quadrances }$q_{1}$\textit{, }$q_{2}%
$\textit{\ and }$q_{3},$\textit{\ corresponding spreads }$S_{1},S_{2}%
$\textit{\ and }$S_{3}$\textit{, and that }$d$\textit{\ is a point lying on
the line }$a_{1}a_{2},$\textit{\ distinct from }$a_{1}$\textit{\ and }$a_{2}%
$\textit{. Define the quadrances }$r_{1}\equiv q\left(  a_{1},d\right)
$\textit{\ and }$r_{2}\equiv q\left(  a_{2},d\right)  $\textit{, and the
spreads }$R_{1}\equiv S\left(  a_{3}a_{1},a_{3}d\right)  $\textit{\ and
}$R_{2}\equiv S\left(  a_{3}a_{2},a_{3}d\right)  $\textit{. Then}%
\[
\frac{R_{1}}{R_{2}}=\frac{S_{1}}{S_{2}}\frac{r_{1}}{r_{2}}=\frac{q_{1}}{q_{2}%
}\frac{r_{1}}{r_{2}}.
\]

\end{theorem}

%

%TCIMACRO{\FRAME{fhFU}{1.8472in}{1.6604in}{0pt}{\Qcb{Triangle proportions:
%$R_{1}/R_{2}=\left(  S_{1}/S_{2}\right)  \times\left(  r_{1}/r_{2}\right)  $}%
%}{\Qlb{TriangleProportions}}{triangleproportions2.EPS}%
%width 1.8472in;  height 1.6604in;  depth 0pt;  original-width 1.9815in;
%original-height 1.7775in;  cropleft "0";  croptop "1";  cropright "1";
%cropbottom "0";  filename 'TriangleProportions2.EPS';file-properties "XNPEU";}%
%} }%
%BeginExpansion
\begin{figure}[h]%
\centering
\includegraphics[
height=1.6604in,
width=1.8472in
]%
{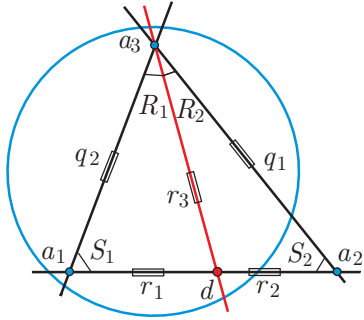}%
\caption{Triangle proportions: $R_{1}/R_{2}=\left(  S_{1}/S_{2}\right)
\times\left(  r_{1}/r_{2}\right)  $}%
\label{TriangleProportions}%
\end{figure}
%EndExpansion

\begin{theorem}
[Menelaus]\textit{Suppose that }$\overline{a_{1}a_{2}a_{3}}$\textit{\ is a
non-null triangle, and that }$L$\textit{\ is a non-null line meeting }%
$a_{2}a_{3}$\textit{, }$a_{1}a_{3}$\textit{\ and }$a_{1}a_{2}$\textit{\ at the
points }$b_{1}$\textit{, }$b_{2}$\textit{\ and }$b_{3}$\textit{\ respectively.
Define the quadrances}%
\[%
\begin{array}
[c]{ccc}%
r_{1}\equiv q\left(  a_{2},b_{1}\right)  &  & t_{1}\equiv q\left(  b_{1}%
,a_{3}\right) \\
r_{2}\equiv q\left(  a_{3},b_{2}\right)  &  & t_{2}\equiv q\left(  b_{2}%
,a_{1}\right) \\
r_{3}\equiv q\left(  a_{1},b_{3}\right)  &  & t_{3}\equiv q\left(  b_{3}%
,a_{2}\right)  .
\end{array}
\]
\textit{Then }$r_{1}r_{2}r_{3}=t_{1}t_{2}t_{3}.$%
%TCIMACRO{\FRAME{fhFU}{5.8033cm}{4.7721cm}{0pt}{\Qcb{Menelaus' theorem:
%$r_{1}r_{2}r_{3}=t_{1}t_{2}t_{3}$}}{}{menelaus.eps}%
%width 5.8033cm;  height 4.7721cm;  depth 0pt;  original-width 6.3668cm;
%original-height 5.3079cm;  cropleft "0";  croptop "1";  cropright "1";
%cropbottom "0";  filename 'Menelaus.eps';file-properties "XNPEU";}} }%
%BeginExpansion
\begin{figure}[h]%
\centering
\includegraphics[
height=4.7721cm,
width=5.8033cm
]%
{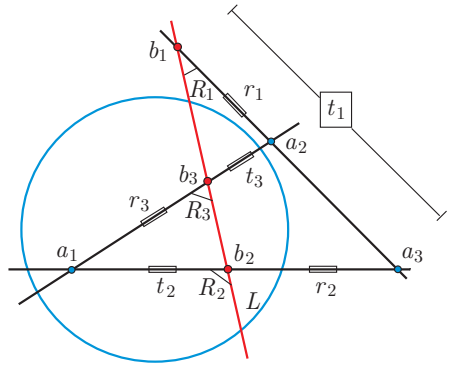}%
\caption{Menelaus' theorem: $r_{1}r_{2}r_{3}=t_{1}t_{2}t_{3}$}%
\end{figure}
%EndExpansion

\end{theorem}

\begin{theorem}
[Menelaus' dual]\textit{Suppose that }$\overline{A_{1}A_{2}A_{3}}$\textit{\ is
a non-null trilateral, and that }$l$\textit{\ is a non-null point joining
}$A_{2}A_{3}$\textit{, }$A_{1}A_{3}$\textit{\ and }$A_{1}A_{2}$\textit{\ on
the lines }$B_{1}$\textit{, }$B_{2}$\textit{\ and }$B_{3}$%
\textit{\ respectively. Define the spreads}%
\[%
\begin{array}
[c]{ccc}%
R_{1}\equiv S\left(  A_{2},B_{1}\right)  &  & T_{1}\equiv S\left(  B_{1}%
,A_{3}\right) \\
R_{2}\equiv S\left(  A_{3},B_{2}\right)  &  & T_{2}\equiv S\left(  B_{2}%
,A_{1}\right) \\
R_{3}\equiv S\left(  A_{1},B_{3}\right)  &  & T_{3}\equiv S\left(  B_{3}%
,A_{2}\right)  .
\end{array}
\]
\textit{Then }$R_{1}R_{2}R_{3}=T_{1}T_{2}T_{3}.$%
%TCIMACRO{\FRAME{fhFU}{4.5612cm}{4.3187cm}{0pt}{\Qcb{Menelaus dual theorem:
%$R_{1}R_{2}R_{3}=T_{1}T_{2}T_{3}$}}{\Qlb{MenelausDual0}}{menelausdual.eps}%
%width 4.5612cm;  height 4.3187cm;  depth 0pt;  original-width 6.7602cm;
%original-height 6.4001cm;  cropleft "0";  croptop "1";  cropright "1";
%cropbottom "0";  filename 'MenelausDual.eps';file-properties "XNPEU";}} }%
%BeginExpansion
\begin{figure}[h]%
\centering
\includegraphics[
height=4.3187cm,
width=4.5612cm
]%
{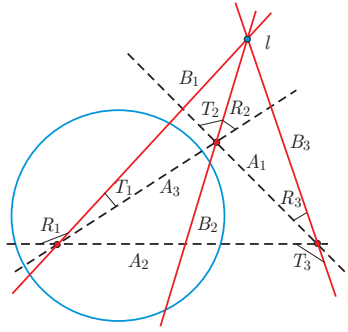}%
\caption{Menelaus dual theorem: $R_{1}R_{2}R_{3}=T_{1}T_{2}T_{3}$}%
\label{MenelausDual0}%
\end{figure}
%EndExpansion

\end{theorem}

\begin{theorem}
[Ceva]\textit{Suppose that the triangle }$\overline{a_{1}a_{2}a_{3}}%
$\textit{\ has non-null lines, that }$a_{0}$\textit{\ is a point distinct from
}$a_{1},a_{2}$\textit{\ and }$a_{3}$\textit{, and that the lines }$a_{0}a_{1}%
$\textit{, }$a_{0}a_{2}$\textit{\ and }$a_{0}a_{3}$\textit{\ meet the lines
}$a_{2}a_{3}$\textit{, }$a_{1}a_{3}$\textit{\ and }$a_{1}a_{2}$%
\textit{\ respectively at the points }$b_{1}$\textit{, }$b_{2}$\textit{\ and
}$b_{3}$\textit{. Define the quadrances}%
\[%
\begin{array}
[c]{ccc}%
r_{1}\equiv q\left(  a_{2},b_{1}\right)  &  & t_{1}\equiv q\left(  b_{1}%
,a_{3}\right) \\
r_{2}\equiv q\left(  a_{3},b_{2}\right)  &  & t_{2}\equiv q\left(  b_{2}%
,a_{1}\right) \\
r_{3}\equiv q\left(  a_{1},b_{3}\right)  &  & t_{3}\equiv q\left(  b_{3}%
,a_{2}\right)  .
\end{array}
\]
\textit{Then }$r_{1}r_{2}r_{3}=t_{1}t_{2}t_{3}.$%
%TCIMACRO{\FRAME{fhFU}{1.8331in}{1.4164in}{0pt}{\Qcb{Ceva's theorem:
%$r_{1}r_{2}r_{3}=t_{1}t_{2}t_{3}$}}{\Qlb{Ceva}}{ceva2.EPS}%
%width 1.8331in;  height 1.4164in;  depth 0pt;  original-width 2.3062in;
%original-height 1.7757in;  cropleft "0";  croptop "1";  cropright "1";
%cropbottom "0";  filename 'Ceva2.EPS';file-properties "XNPEU";}} }%
%BeginExpansion
\begin{figure}[h]%
\centering
\includegraphics[
height=1.4164in,
width=1.8331in
]%
{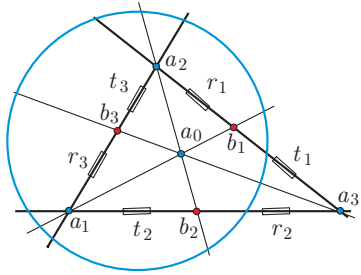}%
\caption{Ceva's theorem: $r_{1}r_{2}r_{3}=t_{1}t_{2}t_{3}$}%
\label{Ceva}%
\end{figure}
%EndExpansion

\end{theorem}

\begin{theorem}
[Ceva dual]\textit{Suppose that the trilateral }$\overline{A_{1}A_{2}A_{3}}%
$\textit{\ is non-null, and that }$A_{0}$\textit{\ is a line distinct from
}$A_{1},A_{2}$\textit{\ and }$A_{3}$\textit{, and that the points }$A_{0}%
A_{1}$\textit{, }$A_{0}A_{2}$\textit{\ and }$A_{0}A_{3}$\textit{\ join the
points }$A_{2}A_{3}$\textit{, }$A_{1}A_{3}$\textit{\ and }$A_{1}A_{2}%
$\textit{\ respectively on the lines }$B_{1}$\textit{, }$B_{2}$\textit{\ and
}$B_{3}$\textit{. Define the spreads}%
\[%
\begin{array}
[c]{ccc}%
R_{1}\equiv S\left(  A_{2},B_{1}\right)  &  & T_{1}\equiv S\left(  B_{1}%
,A_{3}\right) \\
R_{2}\equiv S\left(  A_{3},B_{2}\right)  &  & T_{2}\equiv S\left(  B_{2}%
,A_{1}\right) \\
R_{3}\equiv S\left(  A_{1},B_{3}\right)  &  & T_{3}\equiv S\left(  B_{3}%
,A_{2}\right)  .
\end{array}
\]
\textit{Then }$R_{1}R_{2}R_{3}=T_{1}T_{2}T_{3}.$%
%TCIMACRO{\FRAME{fhFU}{4.7152cm}{4.0509cm}{0pt}{\Qcb{Ceva's dual theorem:
%$R_{1}R_{2}R_{3}=T_{1}T_{2}T_{3}$}}{}{CevaDual.eps}%
%width 4.7152cm;  height 4.0509cm;  depth 0pt;  original-width 6.4825cm;
%original-height 5.5279cm;  cropleft "0";  croptop "1";  cropright "1";
%cropbottom "0";  filename 'CevaDual.eps';file-properties "XNPEU";}} }%
%BeginExpansion
\begin{figure}[h]%
\centering
\includegraphics[
height=4.0509cm,
width=4.7152cm
]%
{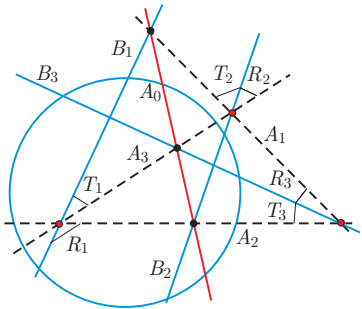}%
\caption{Ceva's dual theorem: $R_{1}R_{2}R_{3}=T_{1}T_{2}T_{3}$}%
\end{figure}
%EndExpansion

\end{theorem}

\section{Isosceles triangles}

\begin{theorem}
[Pons Asinorum]\textit{Suppose that the non-null triangle }$\overline
{a_{1}a_{2}a_{3}}$\textit{\ has quadrances }$q_{1},q_{2}$\textit{\ and }%
$q_{3},$\textit{\ and corresponding spreads }$S_{1},S_{2}$\textit{\ and
}$S_{3}.$ Then\textit{\ }$q_{1}=q_{2}$\textit{\ precisely when }$S_{1}=S_{2}.$
\end{theorem}

\begin{theorem}
[Isosceles right]\textit{If }$\overline{a_{1}a_{2}a_{3}}$\textit{\ is an
isosceles triangle with two right spreads }$S_{1}=S_{2}=1,$\textit{\ then also
}$q_{1}=q_{2}=1$\textit{\ and }$S_{3}=q_{3}.$%
%TCIMACRO{\FRAME{fhFU}{4.2043cm}{3.1917cm}{0pt}{\Qcb{Isosceles right
%triangle:\ $q_{1}=q_{2}=1$ and $S_{3}=q_{3}$}}{}{pythagorasisos.eps}%
%width 4.2043cm;  height 3.1917cm;  depth 0pt;  original-width 6.3908cm;
%original-height 4.8348cm;  cropleft "0";  croptop "1";  cropright "1";
%cropbottom "0";  filename 'PythagorasIsos.eps';file-properties "XNPEU";}} }%
%BeginExpansion
\begin{figure}[h]%
\centering
\includegraphics[
height=3.1917cm,
width=4.2043cm
]%
{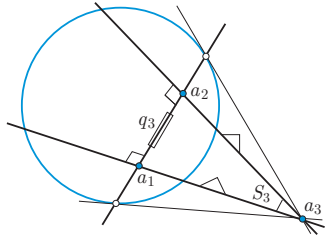}%
\caption{Isosceles right triangle:\ $q_{1}=q_{2}=1$ and $S_{3}=q_{3}$}%
\end{figure}
%EndExpansion

\end{theorem}

\begin{theorem}
[Isosceles triangle]\textit{Suppose a non-null isosceles triangle }%
$\overline{a_{1}a_{2}a_{3}}$\textit{\ has quadrances }$q_{1}=q_{2}\equiv
q$\textit{\ and }$q_{3},$\textit{\ and corresponding spreads }$S_{1}%
=S_{2}\equiv S$\textit{\ and }$S_{3}.$\textit{\ Then the following relations
hold:}%
\[
q_{3}=\frac{4\left(  1-S\right)  q\left(  1-q\right)  }{\left(  1-Sq\right)
^{2}}\qquad\mathrm{and}\qquad S_{3}=\frac{4S\left(  1-S\right)  \left(
1-q\right)  }{\left(  1-Sq\right)  ^{2}}.
\]

\end{theorem}

%

%TCIMACRO{\FRAME{fhFU}{1.4396in}{1.3633in}{0pt}{\Qcb{An isosceles triangle:
%$q_{1}=q_{2}=q,$ $S_{1}=S_{2}=S$}}{\Qlb{IsoscelesTriangle2}}%
%{isoscelestriangle2.EPS}{\special{ language "Scientific Word";
%type "GRAPHIC";  maintain-aspect-ratio TRUE;  display "USEDEF";
%valid_file "F";  width 1.4396in;  height 1.3633in;  depth 0pt;
%original-width 2.1492in;  original-height 2.0338in;  cropleft "0";
%filename 'IsoscelesTriangle2.EPS';file-properties "XNPEU";}} }%
%BeginExpansion
\begin{figure}[h]%
\centering
\includegraphics[
height=1.3633in,
width=1.4396in
]%
{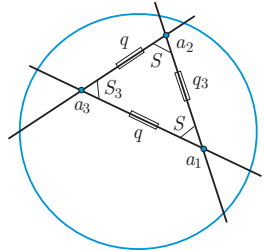}%
\caption{An isosceles triangle: $q_{1}=q_{2}=q,$ $S_{1}=S_{2}=S$}%
\label{IsoscelesTriangle2}%
\end{figure}
%EndExpansion

\begin{theorem}
[Isosceles parallax]\textit{If }$\overline{a_{1}a_{2}a_{3}}$\textit{\ is an
isosceles triangle with }$a_{1}$\textit{\ a null point, }$q_{1}\equiv
q$\textit{\ and }$S_{2}=S_{3}\equiv S,$\textit{\ then }%
\[
q=\frac{4\left(  S-1\right)  }{S^{2}}.
\]

\end{theorem}

%

%TCIMACRO{\FRAME{fhFU}{1.6001in}{1.4094in}{0pt}{\Qcb{Isosceles parallax:
%$q=4\left(  S-1\right)  /S^{2}$}}{\Qlb{IsoscelesParallax}}%
%{isoscelesparallax.EPS}{\special{ language "Scientific Word";
%type "GRAPHIC";  maintain-aspect-ratio TRUE;  display "USEDEF";
%valid_file "F";  width 1.6001in;  height 1.4094in;  depth 0pt;
%original-width 1.9123in;  original-height 1.6799in;  cropleft "0";
%filename 'IsoscelesParallax.EPS';file-properties "XNPEU";}} }%
%BeginExpansion
\begin{figure}[h]%
\centering
\includegraphics[
height=1.4094in,
width=1.6001in
]%
{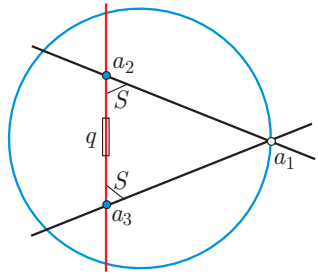}%
\caption{Isosceles parallax: $q=4\left(  S-1\right)  /S^{2}$}%
\label{IsoscelesParallax}%
\end{figure}
%EndExpansion

\section{Equilateral triangles}

\begin{theorem}
[Equilateral quadrance spread]Suppose that a triangle $\overline{a_{1}%
a_{2}a_{3}}$ is equilateral with common quadrance $q_{1}=q_{2}=q_{3}\equiv q,$
and with common spread $S_{1}=S_{2}=S_{3}\equiv S$. Then%
\[
\left(  1-Sq\right)  ^{2}=4\left(  1-S\right)  \left(  1-q\right)  .
\]

\end{theorem}

%

%TCIMACRO{\FRAME{fhFU}{3.7816in}{1.5666in}{0pt}{\Qcb{Equilateral quadrance
%spread theorem: $\left(  1-Sq\right)  ^{2}=4\left(  1-S\right)  \left(
%1-q\right)  $}}{\Qlb{EquilateralBoth}}{equilateralboth2.EPS}%
%width 3.7816in;  height 1.5666in;  depth 0pt;  original-width 4.841in;
%original-height 1.99in;  cropleft "0";  croptop "1";  cropright "1";
%cropbottom "0";  filename 'EquilateralBoth2.EPS';file-properties "XNPEU";}} }%
%BeginExpansion
\begin{figure}[h]%
\centering
\includegraphics[
height=1.5666in,
width=3.7816in
]%
{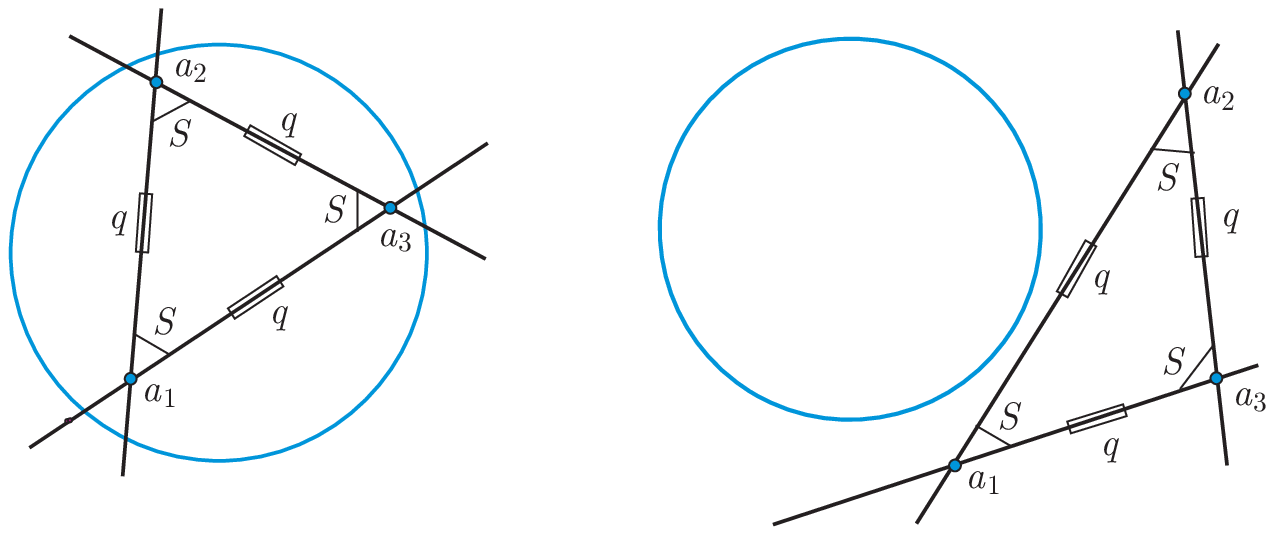}%
\caption{Equilateral quadrance spread theorem: $\left(  1-Sq\right)
^{2}=4\left(  1-S\right)  \left(  1-q\right)  $}%
\label{EquilateralBoth}%
\end{figure}
%EndExpansion

\section{Lambert quadrilaterals}

\begin{theorem}
[Lambert quadrilateral]Suppose a quadrilateral $\overline{abcd}$ has all three
spreads at $a,b$ and $c$ equal to $1.$ Suppose that $q\equiv q\left(
a,b\right)  $ and $p\equiv q\left(  b,c\right)  $. Then%
\begin{align*}
q\left(  c,d\right)   &  =y=\frac{q\left(  1-p\right)  }{1-qp}~\qquad~q\left(
a,d\right)  =x=\frac{p\left(  1-q\right)  }{1-qp}\\
q\left(  a,c\right)   &  =s=q+p-qp~\qquad~q\left(  b,d\right)  =r=\frac
{q+p-2qp}{1-qp}%
\end{align*}
and%
\begin{align*}
S\left(  ba,bd\right)   &  =\frac{x}{r}~\qquad~S\left(  bc,bd\right)
=\frac{y}{r}~\qquad~S\left(  cb,ca\right)  =\frac{q}{s}\\
S\left(  ac,ab\right)   &  =\frac{p}{s}~\qquad~S\left(  ac,ad\right)
=\frac{q\left(  1-p\right)  }{s}~\qquad~S\left(  ca,cd\right)  =\frac{p\left(
1-q\right)  }{s}%
\end{align*}
and%
\[
S\left(  da,dc\right)  =S=1-pq.
\]

\end{theorem}

%

%TCIMACRO{\FRAME{fhFU}{3.9115in}{1.6656in}{0pt}{\Qcb{Lambert quadrilateral
%$\overline{abcd}$}}{\Qlb{Lambert}}{lambert.eps}%
%width 3.9115in;  height 1.6656in;  depth 0pt;  original-width 5.7519in;
%original-height 2.4336in;  cropleft "0";  croptop "1";  cropright "1";
%cropbottom "0";  filename 'Lambert.eps';file-properties "XNPEU";}} }%
%BeginExpansion
\begin{figure}[h]%
\centering
\includegraphics[
height=1.6656in,
width=3.9115in
]%
{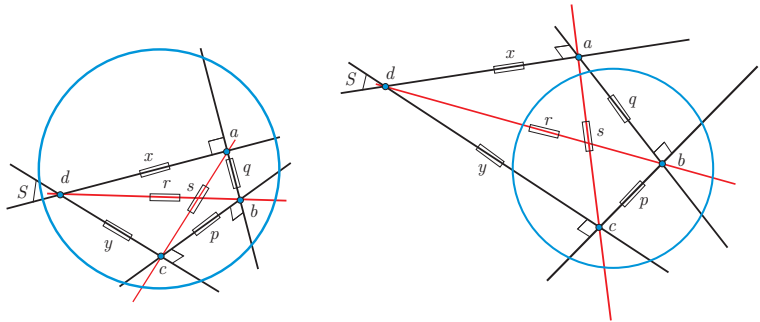}%
\caption{Lambert quadrilateral $\overline{abcd}$}%
\label{Lambert}%
\end{figure}
%EndExpansion

\section{Quadrea and triangle thinness}

\begin{figure}[h]%
\centering
\includegraphics[
height=1.2486in,
width=1.966in
]%
{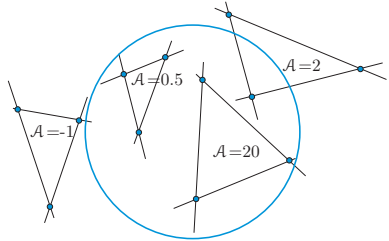}%
\caption{Examples of triangles with quadreas $\mathcal{A}=-1,0.5,20$ and $2$}%
\label{QuadreaExamples}%
\end{figure}

If $\overline{a_{1}a_{2}a_{3}}$ is a triangle with quadrances $q_{1},q_{2}$
and $q_{3},$ and spreads $S_{1},S_{2}$ and $S_{3},$ then from the Spread law
the quantity
\[
\mathcal{A}\equiv S_{1}q_{2}q_{3}=S_{2}q_{1}q_{3}=S_{3}q_{1}q_{2}%
\]
is well-defined, and called the \textbf{quadrea }of the triangle
$\overline{a_{1}a_{2}a_{3}}.$ It is the analog of the squared area in
universal hyperbolic geometry. In Figure \ref{QuadreaExamples} several
triangles with their associated quadreas are shown. Note that the quadrea is
positive for a triangle of internal points, but may also be negative
otherwise.%
%TCIMACRO{\FRAME{fhFU}{1.966in}{1.2486in}{0pt}{\Qcb{Examples of triangles with
%quadreas $\QTR{cal}{A}=-1,0.5,20$ and $2$}}{\Qlb{QuadreaExamples}%
%}{quadreaexamples.EPS}{\special{ language "Scientific Word";  type "GRAPHIC";
%width 1.966in;  height 1.2486in;  depth 0pt;  original-width 3.2021in;
%original-height 2.0232in;  cropleft "0";  croptop "1";  cropright "1";
%cropbottom "0";  filename 'QuadreaExamples.EPS';file-properties "XNPEU";}} }%
%BeginExpansion
%EndExpansion

An interesting aspect of hyperbolic geometry is that \textit{triangles are
thin}. Here are two ways of giving meaning to this, both involving the quadrea
of a triangle.

\begin{theorem}
[Triply nil Cevian thinness]Suppose that $\overline{\alpha_{1}\alpha_{2}%
\alpha_{3}}$ is a triply nil triangle, and that $a$ is a point distinct from
$\alpha_{1},\alpha_{2}$ and $\alpha_{3}.$ Define the cevian points
$c_{1}\equiv\left(  a\alpha_{1}\right)  \left(  \alpha_{2}\alpha_{3}\right)
$, $c_{2}\equiv\left(  a\alpha_{2}\right)  \left(  \alpha_{1}\alpha
_{3}\right)  $ and $c_{3}\equiv\left(  a\alpha_{3}\right)  \left(  \alpha
_{1}\alpha_{2}\right)  .$ Then $\mathcal{A}\left(  \overline{c_{1}c_{2}c_{3}%
}\right)  =1.$%
%TCIMACRO{\FRAME{fhFU}{3.5612in}{1.4014in}{0pt}{\Qcb{Cevian triangle thinness:
%$\QTR{cal}{A}\left(  \overline{c_{1}c_{2}c_{3}}\right)  =1$}}%
%{\Qlb{CevianThinness}}{quadreacevian.EPS}%
%width 3.5612in;  height 1.4014in;  depth 0pt;  original-width 5.4133in;
%original-height 2.1128in;  cropleft "0";  croptop "1";  cropright "1";
%cropbottom "0";  filename 'QuadreaCevian.EPS';file-properties "XNPEU";}} }%
%BeginExpansion
\begin{figure}[h]%
\centering
\includegraphics[
height=1.4014in,
width=3.5612in
]%
{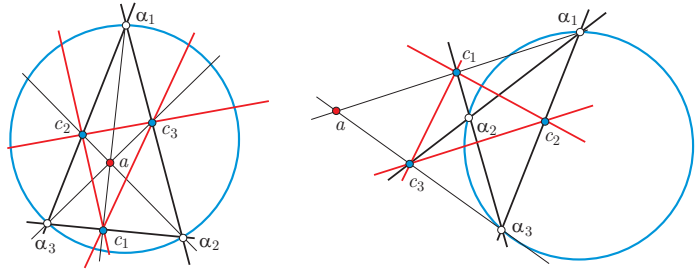}%
\caption{Cevian triangle thinness: $\mathcal{A}\left(  \overline{c_{1}%
c_{2}c_{3}}\right)  =1$}%
\label{CevianThinness}%
\end{figure}
%EndExpansion

\end{theorem}

\begin{theorem}
[Triply nil altitude thinness]Suppose that $\overline{\alpha_{1}\alpha
_{2}\alpha_{3}}$ is a triply nil triangle and that $a$ is a point distinct
from the duals of the lines. If the altitudes to the lines of this triangle
from $a$ meet the lines respectively at base points $b_{1},b_{2}$ and $b_{3}$,
then $\mathcal{A}\left(  \overline{b_{1}b_{2}b_{3}}\right)  =1.$%
%TCIMACRO{\FRAME{fhFU}{3.0787in}{1.4263in}{0pt}{\Qcb{Altitude triangle
%thinness: $\QTR{cal}{A}\left(  \overline{b_{1}b_{2}b_{3}}\right)  =1$}}%
%{}{quadreaaltitudethinness.EPS}{\special{ language "Scientific Word";
%type "GRAPHIC";  maintain-aspect-ratio TRUE;  display "USEDEF";
%valid_file "F";  width 3.0787in;  height 1.4263in;  depth 0pt;
%original-width 4.4434in;  original-height 2.0453in;  cropleft "0";
%filename 'QuadreaAltitudeThinness.EPS';file-properties "XNPEU";}} }%
%BeginExpansion
\begin{figure}[h]%
\centering
\includegraphics[
height=1.4263in,
width=3.0787in
]%
{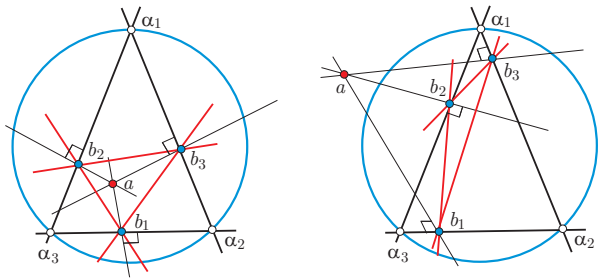}%
\caption{Altitude triangle thinness: $\mathcal{A}\left(  \overline{b_{1}%
b_{2}b_{3}}\right)  =1$}%
\end{figure}
%EndExpansion

\end{theorem}

\section{Null perspective and null subtended theorems}

There are many trigonometric results that prominently feature \textit{null
points }and \textit{null lines}. We give a sample of these now. For some we
include the dual formulations, for others these are left to the reader.

\begin{theorem}
[Null perspective]Suppose that $\alpha_{1},\alpha_{2}$ and $\alpha_{3}$ are
distinct null points, and $b$ is any point on $\alpha_{1}\alpha_{3}$ distinct
from $\alpha_{1}$ and $\alpha_{3}.$ Suppose further that $x$ and $y$ are
points lying on $\alpha_{1}\alpha_{2},$ and that $x_{1}\equiv\left(
\alpha_{2}\alpha_{3}\right)  \left(  xb\right)  $ and $y_{1}\equiv\left(
\alpha_{2}\alpha_{3}\right)  \left(  yb\right)  $. Then%
\[
q\left(  x,y\right)  =q\left(  x_{1},y_{1}\right)  .
\]

\end{theorem}

%

%TCIMACRO{\FRAME{fhFU}{5.4578in}{2.1223in}{0pt}{\Qcb{Null perspective theorem:
%$q\left(  x,y\right)  =q\left(  x_{1},y_{1}\right)  $}}{\Qlb{Null Perspective}%
%}{nullperspectiveboth.EPS}{\special{ language "Scientific Word";
%type "GRAPHIC";  maintain-aspect-ratio TRUE;  display "USEDEF";
%valid_file "F";  width 5.4578in;  height 2.1223in;  depth 0pt;
%original-width 5.8565in;  original-height 2.261in;  cropleft "0";
%filename 'NullPerspectiveBoth.EPS';file-properties "XNPEU";}} }%
%BeginExpansion
\begin{figure}[h]%
\centering
\includegraphics[
height=2.1223in,
width=5.4578in
]%
{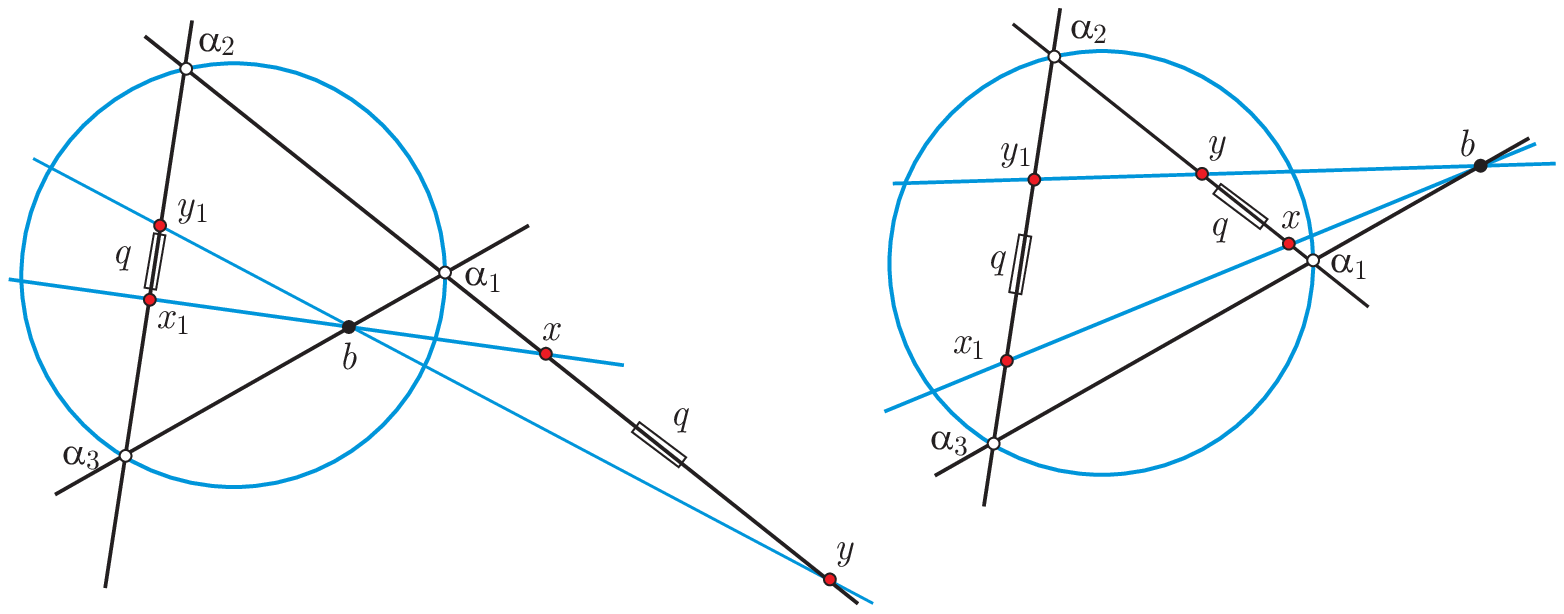}%
\caption{Null perspective theorem: $q\left(  x,y\right)  =q\left(  x_{1}%
,y_{1}\right)  $}%
\label{Null Perspective}%
\end{figure}
%EndExpansion

\begin{theorem}
[Null subtended]\textit{Suppose that the line }$L$\textit{\ passes through the
null points }$\alpha_{1}$\textit{\ and }$\alpha_{2}.$\textit{\ Then for any
other null point }$\alpha_{3}$\textit{\ and any line }$M$\textit{, let }%
$a_{1}\equiv\left(  \alpha_{1}\alpha_{3}\right)  M$\textit{\ and }$a_{2}%
\equiv\left(  \alpha_{2}\alpha_{3}\right)  M.$\textit{\ Then }$q\equiv
q\left(  a_{1},a_{2}\right)  $\textit{\ and }$S\equiv S\left(  L,M\right)
$\textit{\ are related by}%
\[
qS=1.
\]
\textit{In particular }$q$\textit{\ is independent of }$\alpha_{3}.$%
%TCIMACRO{\FRAME{fhFU}{11.7212cm}{4.793cm}{0pt}{\Qcb{Null subtended theorem:
%$qS=1$}}{\Qlb{NullSubtended1}}{nullsubtendedboth.eps}%
%width 11.7212cm;  height 4.793cm;  depth 0pt;  original-width 11.2858cm;
%original-height 4.5712cm;  cropleft "0";  croptop "1";  cropright "1";
%cropbottom "0";  filename 'NullSubtendedBoth.eps';file-properties "XNPEU";}}
%}%
%BeginExpansion
\begin{figure}[h]%
\centering
\includegraphics[
height=4.793cm,
width=11.7212cm
]%
{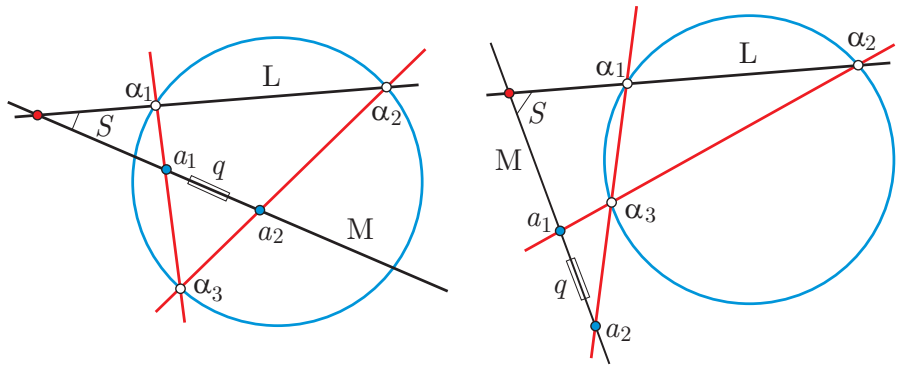}%
\caption{Null subtended theorem: $qS=1$}%
\label{NullSubtended1}%
\end{figure}
%EndExpansion

\end{theorem}

Figure \ref{NullSubtended1} shows two different examples; note that $M$ need
not pass through any null points. The Null subtended theorem allows you to
create a hyperbolic ruler using just a straight-edge, in the sense that you
can use it to repeatedly duplicate a given segment on a given line. Here is
the dual result.

\begin{theorem}
[Null subtended dual]\textit{Suppose that the point }$l$\textit{\ lies on the
null lines }$\Lambda_{1}$\textit{\ and }$\Lambda_{2}.$\textit{\ Then for any
other null line }$\Lambda_{3}$\textit{\ and any point }$m$\textit{, let
}$A_{1}\equiv\left(  \Lambda_{1}\Lambda_{3}\right)  m$\textit{\ and }%
$A_{2}\equiv\left(  \Lambda_{2}\Lambda_{3}\right)  m.$\textit{\ Then }$S\equiv
S\left(  L_{1},L_{2}\right)  $\textit{\ and }$q\equiv q\left(  l,m\right)
$\textit{\ are related by}%
\[
Sq=1.
\]

\end{theorem}

\begin{theorem}
[Opposite subtended]\textit{Suppose }$\overline{\alpha\beta\gamma\delta}%
$\textit{\ is a quadrangle of null points, and that }$\upsilon,\mu
$\textit{\ are also null points. Let }$a\equiv\left(  \alpha\mu\right)
\left(  \gamma\delta\right)  $\textit{, }$b\equiv\left(  \beta\mu\right)
\left(  \gamma\delta\right)  $\textit{, }$c\equiv\left(  \gamma\upsilon
\right)  \left(  \alpha\beta\right)  $\textit{\ and }$d\equiv\left(
\delta\upsilon\right)  \left(  \alpha\beta\right)  .$\textit{\ Then}%
\[
q\left(  a,b\right)  =q\left(  c,d\right)  .
\]

\end{theorem}

%

%TCIMACRO{\FRAME{fhFU}{5.4982cm}{3.8441cm}{0pt}{\Qcb{Opposite subtended spreads
%theorem: $q\left(  a,b\right)  =q\left(  c,d\right)  $}}%
%{\Qlb{Oppositesubtended1}}{nullsubtendedcor.eps}%
%width 5.4982cm;  height 3.8441cm;  depth 0pt;  original-width 6.092cm;
%original-height 4.2305cm;  cropleft "0";  croptop "1";  cropright "1";
%cropbottom "0";  filename 'NullSubtendedCor.eps';file-properties "XNPEU";}} }%
%BeginExpansion
\begin{figure}[h]%
\centering
\includegraphics[
height=3.8441cm,
width=5.4982cm
]%
{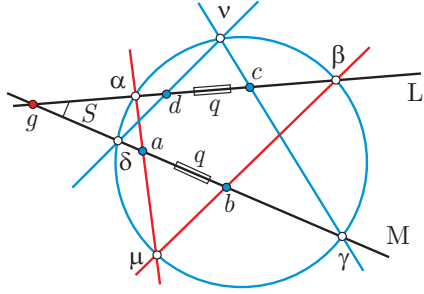}%
\caption{Opposite subtended spreads theorem: $q\left(  a,b\right)  =q\left(
c,d\right)  $}%
\label{Oppositesubtended1}%
\end{figure}
%EndExpansion

Butterfly theorems have been investigated in the hyperbolic plane (\cite{S
J}). The next theorems concern a related configuration of null points.

\begin{theorem}
[Butterfly quadrance]\textit{Suppose that }$\overline{\alpha\beta\gamma\delta
}$\textit{\ is a quadrangle of null points, with }$g\equiv\left(  \alpha
\gamma\right)  \left(  \beta\delta\right)  $\textit{\ a diagonal point. Let
}$L$\textit{\ be any line passing through }$g,$\textit{\ and suppose that }%
$L$\textit{\ meets }$\alpha\delta$\textit{\ at }$x$\textit{\ and }$\beta
\gamma$\textit{\ at }$y$\textit{. Then }%
\[
q\left(  g,x\right)  =q\left(  g,y\right)  .
\]

\end{theorem}

%

%TCIMACRO{\FRAME{fhFU}{3.3546in}{1.3757in}{0pt}{\Qcb{Butterfly quadrance
%theorem: $q\left(  g,x\right)  =q\left(  g,y\right)  $}}%
%{\Qlb{ButterflyQuadrance}}{butterflyboth2.EPS}%
%width 3.3546in;  height 1.3757in;  depth 0pt;  original-width 4.5981in;
%original-height 1.8689in;  cropleft "0";  croptop "1";  cropright "1";
%cropbottom "0";  filename 'ButterflyBoth2.EPS';file-properties "XNPEU";}} }%
%BeginExpansion
\begin{figure}[h]%
\centering
\includegraphics[
height=1.3757in,
width=3.3546in
]%
{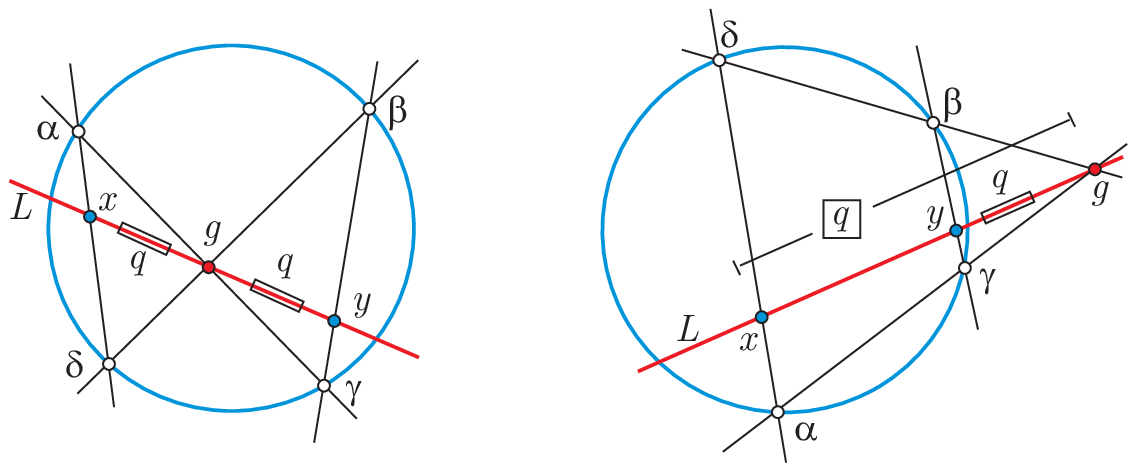}%
\caption{Butterfly quadrance theorem: $q\left(  g,x\right)  =q\left(
g,y\right)  $}%
\label{ButterflyQuadrance}%
\end{figure}
%EndExpansion

\begin{theorem}
[Butterfly spread]\textit{Suppose that }$\overline{\alpha\beta\gamma\delta}%
$\textit{\ is a quadrangle of null points, with }$g\equiv\left(  \alpha
\gamma\right)  \left(  \beta\delta\right)  $\textit{\ a diagonal point. Let
}$L$\textit{\ be any line passing through }$g$\textit{. Then }%
\[
S\left(  L,\alpha\delta\right)  =S\left(  L,\beta\gamma\right)  .
\]

\end{theorem}

%

%TCIMACRO{\FRAME{fhFU}{3.4113in}{1.4342in}{0pt}{\Qcb{Butterfly spread theorem:
%$S\left(  L,\alpha\delta\right)  =S\left(  L,\beta\gamma\right)  $}%
%}{\Qlb{ButterflySpread}}{butterflyspread2.EPS}%
%width 3.4113in;  height 1.4342in;  depth 0pt;  original-width 4.4855in;
%original-height 1.8689in;  cropleft "0";  croptop "1";  cropright "1";
%cropbottom "0";  filename 'ButterflySpread2.EPS';file-properties "XNPEU";}} }%
%BeginExpansion
\begin{figure}[h]%
\centering
\includegraphics[
height=1.4342in,
width=3.4113in
]%
{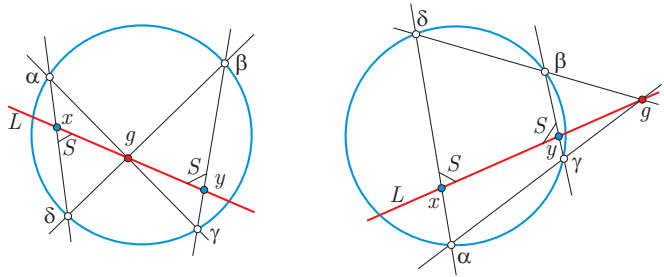}%
\caption{Butterfly spread theorem: $S\left(  L,\alpha\delta\right)  =S\left(
L,\beta\gamma\right)  $}%
\label{ButterflySpread}%
\end{figure}
%EndExpansion

\section{The 48/64 theorems}

In universal hyperbolic geometry we discover many constants of nature that
express themselves in a geometrical way. Prominent among these are the numbers
$48$ and $64$, but there are many others too!

\begin{theorem}
[ $48/64$]\textit{If the three spreads between opposite lines of a quadrangle
}$\overline{\alpha_{1}\alpha_{2}\alpha_{3}\alpha_{4}}$\textit{\ of null points
are }$P,R$\textit{\ and }$T,$\textit{\ then}%
\[
PR+RT+PT=48
\]
and
\[
PRT=64.
\]

\end{theorem}

%

%TCIMACRO{\FRAME{fhFU}{6.9194cm}{5.2346cm}{0pt}{\Qcb{The $48/64$
%theorem:\ $PR+RT+PT=48 $ and $PRT=64$}}{\Qlb{48Theorem1}}%
%{nullquadrangle48.eps}{\special{ language "Scientific Word";  type "GRAPHIC";
%width 6.9194cm;  height 5.2346cm;  depth 0pt;  original-width 6.7971cm;
%original-height 5.1187cm;  cropleft "0";  croptop "1";  cropright "1";
%cropbottom "0";  filename 'NullQuadrangle48.eps';file-properties "XNPEU";}} }%
%BeginExpansion
\begin{figure}[h]%
\centering
\includegraphics[
height=5.2346cm,
width=6.9194cm
]%
{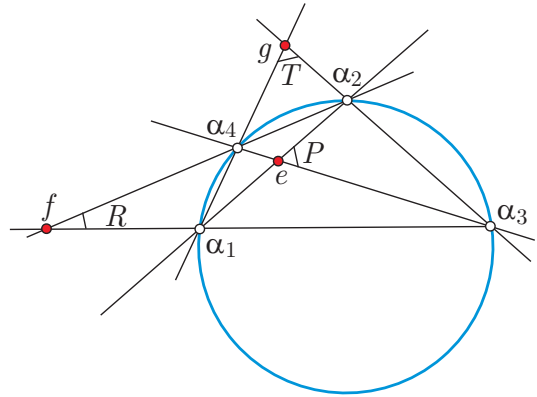}%
\caption{The $48/64$ theorem:\ $PR+RT+PT=48 $ and $PRT=64$}%
\label{48Theorem1}%
\end{figure}
%EndExpansion

It follows that
\[
\frac{1}{R}+\frac{1}{S}+\frac{1}{T}=\frac{3}{4}.
\]
In particular if we know two of these spreads, we get a linear equation for
the third one.

\begin{theorem}
[ $48/64$ dual]\textit{If the three quadrances between opposite points of a
quadrilateral }$\overline{\Lambda_{1}\Lambda_{2}\Lambda_{3}\Lambda_{4}}%
$\textit{\ of null lines are }$p,r$\textit{\ and }$t,$\textit{\ then}%
\[
pr+rt+pt=48
\]
and
\[
prt=64.
\]

\end{theorem}

%

%TCIMACRO{\FRAME{fhFU}{6.2758cm}{5.0083cm}{0pt}{\Qcb{The $48/64$ dual theorem:
%$pr+rt+pt=48$ and $prt=64$}}{\Qlb{48Dual}}{nullquadrangle48dual.eps}%
%width 6.2758cm;  height 5.0083cm;  depth 0pt;  original-width 5.6407cm;
%original-height 4.4844cm;  cropleft "0";  croptop "1";  cropright "1";
%cropbottom "0";  filename 'NullQuadrangle48Dual.eps';file-properties "XNPEU";}%
%} }%
%BeginExpansion
\begin{figure}[h]%
\centering
\includegraphics[
height=5.0083cm,
width=6.2758cm
]%
{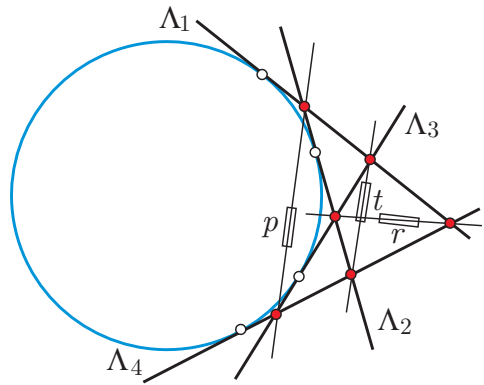}%
\caption{The $48/64$ dual theorem: $pr+rt+pt=48$ and $prt=64$}%
\label{48Dual}%
\end{figure}
%EndExpansion

\section{Pentagon theorems and extensions}

The next theorem does not rely on null points, but is closely connected to a
family of results that do.

\begin{theorem}
[Pentagon ratio]\textit{Suppose }$\overline{a_{1}a_{2}a_{3}a_{4}a_{5}}%
$\textit{\ is a pentagon, meaning a cyclical list of five points, no three
consecutive points collinear. Define }\textbf{diagonal points }$b_{1}%
\equiv\left(  a_{2}a_{4}\right)  \left(  a_{3}a_{5}\right)  $\textit{, }%
$b_{2}\equiv\left(  a_{3}a_{5}\right)  \left(  a_{4}a_{1}\right)  $\textit{,
}$b_{3}\equiv\left(  a_{4}a_{1}\right)  \left(  a_{5}a_{2}\right)  $\textit{,
}$b_{4}\equiv\left(  a_{5}a_{2}\right)  \left(  a_{1}a_{3}\right)
$\textit{\ and }$b_{5}\equiv\left(  a_{1}a_{3}\right)  \left(  a_{2}%
a_{4}\right)  $\textit{, and subsequently }\textbf{opposite points }%
$c_{1}\equiv\left(  a_{1}b_{1}\right)  \left(  a_{2}a_{5}\right)  $\textit{,
}$c_{2}\equiv\left(  a_{2}b_{2}\right)  \left(  a_{3}a_{1}\right)  $\textit{,
}$c_{3}\equiv\left(  a_{3}b_{3}\right)  \left(  a_{4}a_{2}\right)  $\textit{,
}$c_{4}\equiv\left(  a_{4}b_{4}\right)  \left(  a_{5}a_{3}\right)
$\textit{\ and }$c_{5}\equiv\left(  a_{5}b_{5}\right)  \left(  a_{1}%
a_{4}\right)  $\textit{. Then}%
\[
q\left(  b_{1},c_{4}\right)  q\left(  b_{2},c_{5}\right)  q\left(  b_{3}%
,c_{1}\right)  q\left(  b_{4},c_{2}\right)  q\left(  b_{5},c_{3}\right)
=q\left(  b_{2},c_{4}\right)  q\left(  b_{3},c_{5}\right)  q\left(
b_{4},c_{1}\right)  q\left(  b_{5},c_{2}\right)  q\left(  b_{1},c_{3}\right)
.
\]

\end{theorem}

%

%TCIMACRO{\FRAME{fhFU}{5.4718cm}{5.0808cm}{0pt}{\Qcb{Pentagon ratio theorem}%
%}{\Qlb{PentagonRatio}}{pentagonalternating.eps}%
%width 5.4718cm;  height 5.0808cm;  depth 0pt;  original-width 6.7537cm;
%original-height 6.2634cm;  cropleft "0";  croptop "1";  cropright "1";
%cropbottom "0";  filename 'PentagonAlternating.eps';file-properties "XNPEU";}}
%}%
%BeginExpansion
\begin{figure}[h]%
\centering
\includegraphics[
height=5.0808cm,
width=5.4718cm
]%
{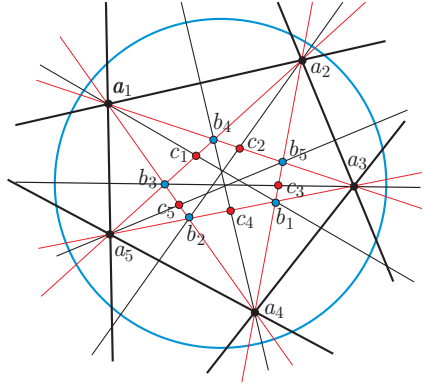}%
\caption{Pentagon ratio theorem}%
\label{PentagonRatio}%
\end{figure}
%EndExpansion

Since the pentagon is arbitrary, it follows by a scaling argument that
\textit{exactly the same theorem} holds in planar Euclidean geometry, where we
replace the hyperbolic quadrance $q$ with the Euclidean quadrance $Q,$ since
for five very close interior points, the hyperbolic quadrances and Euclidean
quadrances are approximately equal.

There are also some interesting additional features that occur in the special
case of the pentagon when all the points $a_{i}$ are null.

\begin{theorem}
[Pentagon null product]\textit{Suppose }$\overline{\alpha_{1}\alpha_{2}%
\alpha_{3}\alpha_{4}\alpha_{5}}$\textit{\ is a pentagon of null points. Define
diagonal points\ }$b_{1}\equiv\left(  \alpha_{2}\alpha_{4}\right)  \left(
\alpha_{3}\alpha_{5}\right)  $\textit{, }$b_{2}\equiv\left(  \alpha_{3}%
\alpha_{5}\right)  \left(  \alpha_{4}\alpha_{1}\right)  $\textit{, }%
$b_{3}\equiv\left(  \alpha_{4}\alpha_{1}\right)  \left(  \alpha_{5}\alpha
_{2}\right)  $\textit{, }$b_{4}\equiv\left(  \alpha_{5}\alpha_{2}\right)
\left(  \alpha_{1}\alpha_{3}\right)  $\textit{\ and }$b_{5}\equiv\left(
\alpha_{1}\alpha_{3}\right)  \left(  \alpha_{2}\alpha_{4}\right)  $\textit{.
Then}%
\[
q\left(  b_{1},b_{2}\right)  q\left(  b_{2},b_{3}\right)  q\left(  b_{3}%
,b_{4}\right)  q\left(  b_{4},b_{5}\right)  q\left(  b_{5},b_{1}\right)
=-\frac{1}{4^{5}}.
\]

\end{theorem}

%

%TCIMACRO{\FRAME{fhFU}{5.3989in}{2.0158in}{0pt}{\Qcb{Pentagon null product
%theorem: $q\left(  b_{1},b_{2}\right)  q\left(  b_{2},b_{3}\right)  q\left(
%b_{3},b_{4}\right)  q\left(  b_{4},b_{5}\right)  q\left(  b_{5},b_{1}\right)
%=-\frac{1}{4^{5}}$}}{\Qlb{Pentagon Null}}{pentagonnullproductboth.EPS}%
%width 5.3989in;  height 2.0158in;  depth 0pt;  original-width 4.7648in;
%original-height 1.7618in;  cropleft "0";  croptop "1";  cropright "1";
%filename 'PentagonNullProductBOTH.EPS';file-properties "XNPEU";}} }%
%BeginExpansion
\begin{figure}[h]%
\centering
\includegraphics[
height=2.0158in,
width=5.3989in
]%
{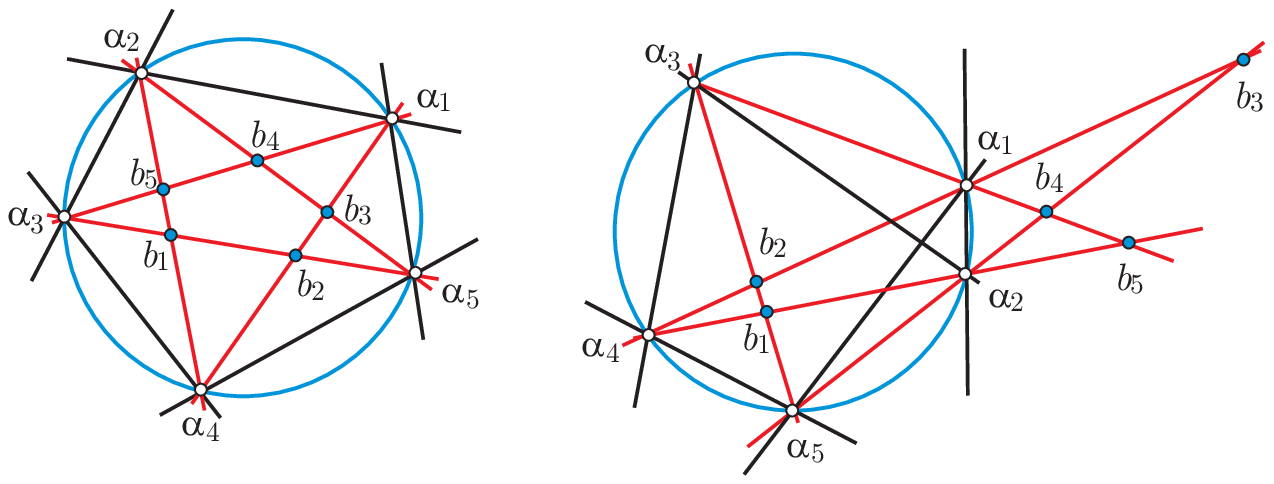}%
\caption{Pentagon null product theorem: $q\left(  b_{1},b_{2}\right)  q\left(
b_{2},b_{3}\right)  q\left(  b_{3},b_{4}\right)  q\left(  b_{4},b_{5}\right)
q\left(  b_{5},b_{1}\right)  =-\frac{1}{4^{5}}$}%
\label{Pentagon Null}%
\end{figure}
%EndExpansion

\begin{theorem}
[Pentagon null symmetry]\textit{With notation as in the Pentagon ratio
theorem, suppose that }$\overline{\alpha_{1}\alpha_{2}\alpha_{3}\alpha
_{4}\alpha_{5}}$\textit{\ is a pentagon of null points, then }%
\begin{align*}
q\left(  b_{1},c_{4}\right)   &  =q\left(  b_{5},c_{2}\right)  \qquad q\left(
b_{2},c_{5}\right)  =q\left(  b_{1},c_{3}\right)  \qquad q\left(  b_{3}%
,c_{1}\right)  =q\left(  b_{2},c_{4}\right) \\
q\left(  b_{4},c_{2}\right)   &  =q\left(  b_{3},c_{5}\right)  \qquad
\mathrm{and}\qquad q\left(  b_{5},c_{3}\right)  =q\left(  b_{4},c_{1}\right)
.
\end{align*}
\textit{Furthermore if we fix }$\alpha_{2},\alpha_{3},\alpha_{4}$\textit{\ and
}$\alpha_{5},$\textit{\ then the quadrance }$q\left(  b_{4},c_{2}\right)
=q\left(  b_{3},c_{5}\right)  $\textit{\ is constant, independent of }%
$\alpha_{1}.$%
%TCIMACRO{\FRAME{fhFU}{2.7273in}{2.0847in}{0pt}{\Qcb{Pentagon null symmetry
%theorem: $q\left(  b_{1},c_{4}\right)  =q\left(  b_{5},c_{2}\right)  $ etc}%
%}{\Qlb{PentagonNullSymmetry}}{pentagonnullsymmetry.EPS}%
%width 2.7273in;  height 2.0847in;  depth 0pt;  original-width 3.563in;
%original-height 2.7168in;  cropleft "0";  croptop "1";  cropright "1";
%cropbottom "0";  filename 'PentagonNullsymmetry.EPS';file-properties "XNPEU";}%
%} }%
%BeginExpansion
\begin{figure}[h]%
\centering
\includegraphics[
height=2.0847in,
width=2.7273in
]%
{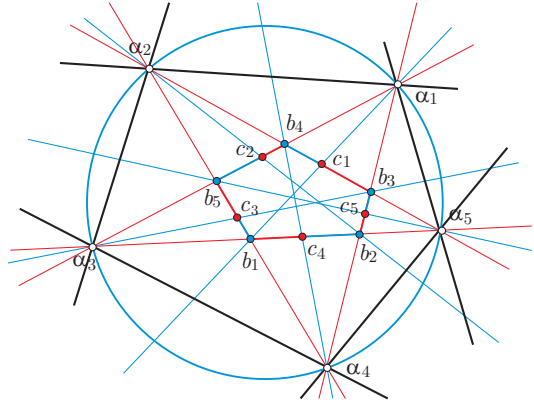}%
\caption{Pentagon null symmetry theorem: $q\left(  b_{1},c_{4}\right)
=q\left(  b_{5},c_{2}\right)  $ etc}%
\label{PentagonNullSymmetry}%
\end{figure}
%EndExpansion

\end{theorem}

Since five points determine a conic, here is an analog to the Pentagon ratio
theorem for general septagons.

\begin{theorem}
[Septagon conic ratio]\textit{Suppose }$\overline{\alpha_{1}\alpha_{2}%
\alpha_{3}\alpha_{4}\alpha_{5}\alpha_{6}\alpha_{7}}$\textit{\ is a septagon of
points lying on a conic. Define diagonal points }$b_{1}\equiv\left(
\alpha_{3}\alpha_{5}\right)  \left(  \alpha_{4}\alpha_{6}\right)  $\textit{,
}$b_{2}\equiv\left(  \alpha_{4}\alpha_{6}\right)  \left(  \alpha_{5}\alpha
_{7}\right)  $\textit{, }$b_{3}\equiv\left(  \alpha_{5}\alpha_{7}\right)
\left(  \alpha_{6}\alpha_{1}\right)  $\textit{, }$b_{4}\equiv\left(
\alpha_{6}\alpha_{1}\right)  \left(  \alpha_{7}\alpha_{2}\right)  $\textit{,
}$b_{5}\equiv\left(  \alpha_{7}\alpha_{2}\right)  \left(  \alpha_{1}\alpha
_{3}\right)  $\textit{, }$b_{6}\equiv\left(  \alpha_{1}\alpha_{3}\right)
\left(  \alpha_{2}\alpha_{4}\right)  $\textit{\ and }$b_{7}\equiv\left(
\alpha_{2}\alpha_{4}\right)  \left(  \alpha_{3}\alpha_{5}\right)  $\textit{,
and opposite points }$c_{1}\equiv\left(  \alpha_{1}b_{1}\right)  \left(
\alpha_{7}\alpha_{2}\right)  $\textit{, }$c_{2}\equiv\left(  \alpha_{2}%
b_{2}\right)  \left(  \alpha_{1}\alpha_{3}\right)  $\textit{, }$c_{3}%
\equiv\left(  \alpha_{3}b_{3}\right)  \left(  \alpha_{2}\alpha_{4}\right)
$\textit{, }$c_{4}\equiv\left(  \alpha_{4}b_{4}\right)  \left(  \alpha
_{3}\alpha_{5}\right)  $\textit{, }$c_{5}\equiv\left(  \alpha_{5}b_{5}\right)
\left(  \alpha_{4}\alpha_{6}\right)  $\textit{, }$c_{6}\equiv\left(
\alpha_{6}b_{6}\right)  \left(  \alpha_{5}\alpha_{7}\right)  $\textit{\ and
}$c_{7}\equiv\left(  \alpha_{7}b_{7}\right)  \left(  \alpha_{6}\alpha
_{1}\right)  $\textit{. Then }%
\begin{align*}
&  q\left(  c_{1},b_{5}\right)  q\left(  c_{2},b_{6}\right)  q\left(
c_{3},b_{7}\right)  q\left(  c_{4},b_{1}\right)  q\left(  c_{5},b_{2}\right)
q\left(  c_{6},b_{3}\right)  q\left(  c_{7},b_{4}\right) \\
&  =q\left(  c_{1},b_{4}\right)  q\left(  c_{2},b_{5}\right)  q\left(
c_{3},b_{6}\right)  q\left(  c_{4},b_{7}\right)  q\left(  c_{5},b_{1}\right)
q\left(  c_{6},b_{2}\right)  q\left(  c_{7},b_{3}\right)  .
\end{align*}

\end{theorem}

%

%TCIMACRO{\FRAME{fhFU}{5.2744in}{2.3138in}{0pt}{\Qcb{Septagon conic ratio
%theorem}}{\Qlb{Heptagon Conic Ratio}}{heptagonconicratio.EPS}%
%width 5.2744in;  height 2.3138in;  depth 0pt;  original-width 5.8164in;
%original-height 2.5389in;  cropleft "0";  croptop "1";  cropright "1";
%cropbottom "0";  filename 'HeptagonConicRatio.EPS';file-properties "XNPEU";}}
%}%
%BeginExpansion
\begin{figure}[h]%
\centering
\includegraphics[
height=2.3138in,
width=5.2744in
]%
{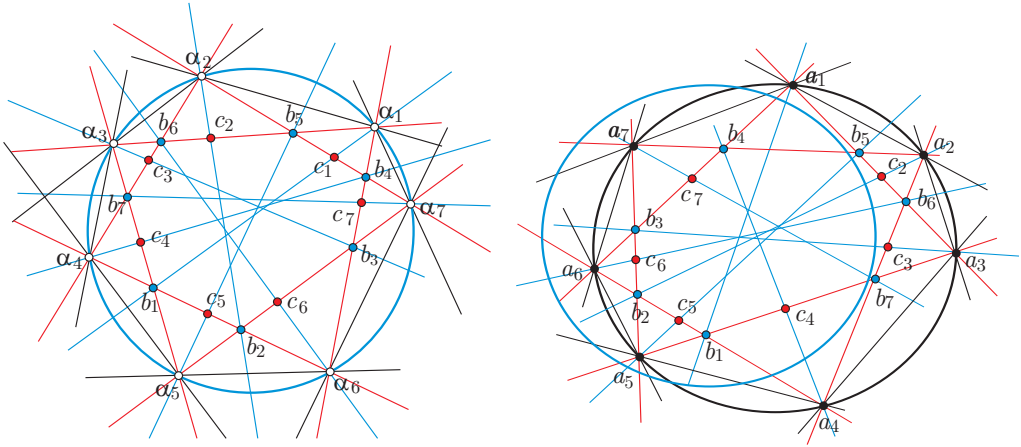}%
\caption{Septagon conic ratio theorem}%
\label{Heptagon Conic Ratio}%
\end{figure}
%EndExpansion

Since the notion of a conic is projective, a scaling argument shows that the
same theorem holds also in the Euclidean case. Figure
\ref{Heptagon Conic Ratio} shows on the left the special case of a septagon of
null points, and on the right a more general case where the septagon lies on a
conic, in this case a Euclidean circle.

I conjecture that the Septagon conic ratio theorem extends to all odd polygons.

\section{Conics in hyperbolic geometry}

The previous result used the fact that conics are well-defined in hyperbolic
geometry, since they can be defined projectively, and we are working in a
projective setting. A natural question is: can we also study conics metrically
as we do in the Euclidean plane? In fact we can, and the resulting theory is
both more intricate and richer than the Euclidean theory, nevertheless
incorporating the Euclidean case as a limiting special case.

\begin{figure}[h]%
\centering
\includegraphics[
height=3.3271in,
width=4.4899in
]%
{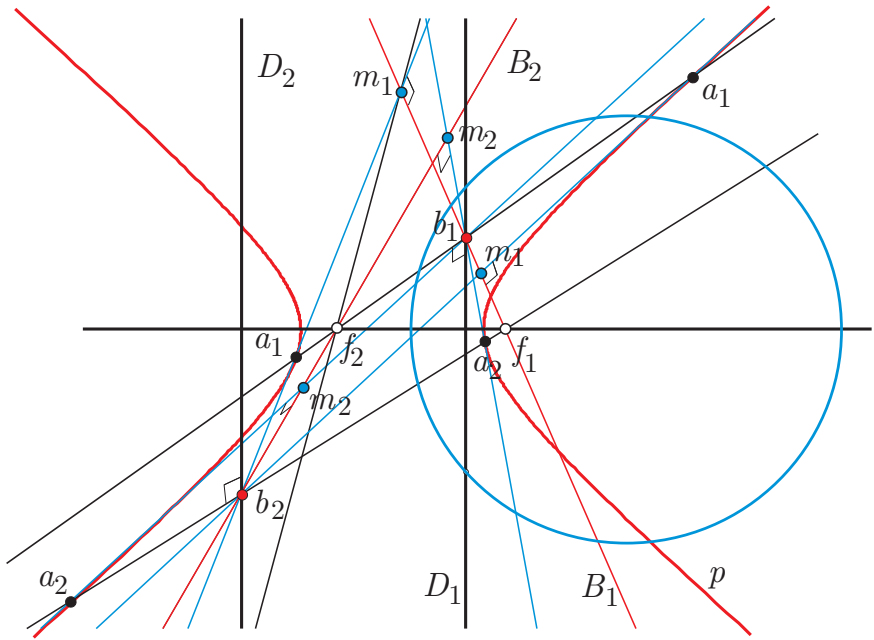}%
\caption{Construction of a (hyperbolic) parabola}%
\label{ParabolaConstruction}%
\end{figure}

We have already mentioned (hyperbolic) circles and illustrated them in Figures
\ref{CirclesInterior} and \ref{CirclesExterior}. Let us now just briefly
outline some results for a\textbf{\ (hyperbolic) parabola},\textbf{\ }which
may be defined as the locus of a point $a$ satisfying $q\left(  a,f\right)
=q\left(  a,D\right)  $ where $f$ is a fixed point called a \textbf{focus},
and $D$ is a fixed line called a \textbf{directrix}, and where $q\left(
a,D\right)  $ is the quadrance from the point $a$ to the base point $b $ of
the altitude to $D$ through $a.$ The following theorems summarize some basic
facts about such a hyperbolic parabola, some similar to the Euclidean
situation, others quite different. The situation is illustrated in Figure
\ref{ParabolaConstruction}. Careful examination reveals many more interesting
features of this situation, which will be discussed in a further paper in this series.

\begin{theorem}
[Parabola focus directrix pair]If a (hyperbolic) parabola $p$ has focus
$f_{1}$ and directrix $D_{1},$ then it also has another focus $f_{2}\equiv
D_{1}^{\perp}$ and another directrix $D_{2}\equiv f_{1}^{\perp}.$
\end{theorem}

\begin{theorem}
[Parabola tangents]Suppose that $b_{1}$ is a point on $D_{1}$ such that the
two midlines of the side $\overline{b_{1}f_{1}}$ exist. Then these midlines
meet the altitude line to $D_{1}$ through $b_{1}$ at two points (both labelled
$a_{1}$ in the Figure) lying on the parabola, and are the tangents to the
parabola at those points.
\end{theorem}

We note that in addition if the two midlines of the side $\overline{b_{1}%
f_{1}}$ exist, then both midlines meet at the point $b_{2}\equiv\left(
b_{1}f_{1}\right)  ^{\perp}$ lying on $D_{2}$ and the corresponding midlines
of the side $\overline{b_{2}f_{2}}$ meet the altitude line to $D_{2}$ through
$b_{2}$ at two points (both labelled $a_{2}$ in the Figure) lying on the
parabola, and themselves meet at $b_{1}.$ This gives a pairing between some of
the points $b_{1}$ lying on $D_{1}$ and some of the points $b_{2}$ lying on
$D_{2}.$ Of the four points labelled $a_{1}$ and $a_{2}$ lying on the
parabola, one of the $a_{1}$ points and one of the $a_{2}$ points are
(somewhat mysteriously) perpendicular. The entire situation is very rich, and
emphasizes once again (see \cite{Wild3}) that the theory of conics is not a
closed book, but rather a rich mine which has only been partly explored so
far.%
%TCIMACRO{\FRAME{fhFU}{4.4899in}{3.3271in}{0pt}{\Qcb{Construction of a
%(hyperbolic) parabola}}{\Qlb{ParabolaConstruction}}{parabolaconstruction2.EPS}%
%width 4.4899in;  height 3.3271in;  depth 0pt;  original-width 4.0003in;
%original-height 2.9554in;  cropleft "0";  croptop "1";  cropright "1";
%filename 'ParabolaConstruction2.EPS';file-properties "XNPEU";}} }%
%BeginExpansion
%EndExpansion

Recall that in Euclidean geometry the locus of a point $a$ satisfying
$q\left(  a,f_{1}\right)  +q\left(  a,f_{2}\right)  =k$ for two fixed points
$f_{1}$ and $f_{2}$ and some fixed number $k$ is a circle.

\begin{theorem}
[Sum of two quadrances]The (hyperbolic) parabola $p$ described in the previous
theorem may also be defined as the locus of those points $a$ satisfying
\[
q\left(  a,f_{1}\right)  +q\left(  a,f_{2}\right)  =1.
\]

\end{theorem}

Many classical theorems for the Euclidean parabola hold also for the
hyperbolic parabola $p.$ Here are two, illustrated in Figure
\ref{Parabola chord}.

\begin{theorem}
[Parabola chord spread]If $a$ and $b$ are two points on the hyperbolic
parabola $p$ with directrix $D$ and focus $f,$ and if $c$ is the meet of $D$
with the tangent of $p$ at $a,$ while $d$ is the meet of $D$ with the tangent
to $p$ at $b,$ then $S\left(  cf,fd\right)  =S\left(  af,fb\right)  $.
\end{theorem}

\begin{theorem}
[Parabola chord tangents perpendicular]If $a$ and $b$ are two points on the
hyperbolic parabola $p$ with directrix $D$ and focus $f$, and if $e$ is the
meet of $D$ with $ab,$ while $g$ is the meet of the tangents to $p$ at $a$ and
$b,$ then $ef$ is perpendicular to $gf.$%
%TCIMACRO{\FRAME{fhFU}{1.7459in}{2.0955in}{0pt}{\Qcb{Hyperbolic parabola with
%focus $f$ and directrix $D$}}{\Qlb{Parabola chord}}{parabolafocus.EPS}%
%width 1.7459in;  height 2.0955in;  depth 0pt;  original-width 3.0592in;
%original-height 3.6783in;  cropleft "0";  croptop "1";  cropright "1";
%cropbottom "0";  filename 'ParabolaFocus.EPS';file-properties "XNPEU";}} }%
%BeginExpansion
\begin{figure}[h]%
\centering
\includegraphics[
height=2.0955in,
width=1.7459in
]%
{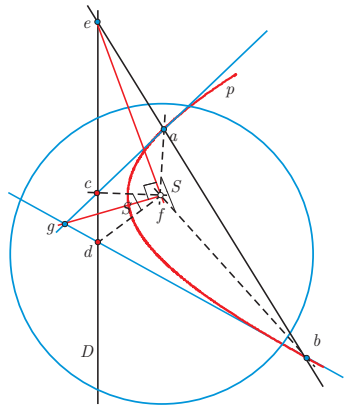}%
\caption{Hyperbolic parabola with focus $f$ and directrix $D$}%
\label{Parabola chord}%
\end{figure}
%EndExpansion

\end{theorem}

\section{Bolyai's construction of limiting lines}

Here is a universal version of a famous construction of J. Bolyai, to find the
limiting lines $U$ and $V$ to an interior line $L\ $through a point $a,$ where
limiting means that $U$ and $V$ meet $L$ on the null circle.%
%TCIMACRO{\FRAME{fhFU}{3.9601in}{2.0357in}{0pt}{\Qcb{A variant on J. Bolyai's
%construction of the limiting lines from $a$ to $L$}}{}{bolyaiparallel2.EPS}%
%width 3.9601in;  height 2.0357in;  depth 0pt;  original-width 3.77in;
%original-height 1.9253in;  cropleft "0";  croptop "1";  cropright "1";
%cropbottom "0";  filename 'BolyaiParallel2.EPS';file-properties "XNPEU";}} }%
%BeginExpansion
\begin{figure}[h]%
\centering
\includegraphics[
height=2.0357in,
width=3.9601in
]%
{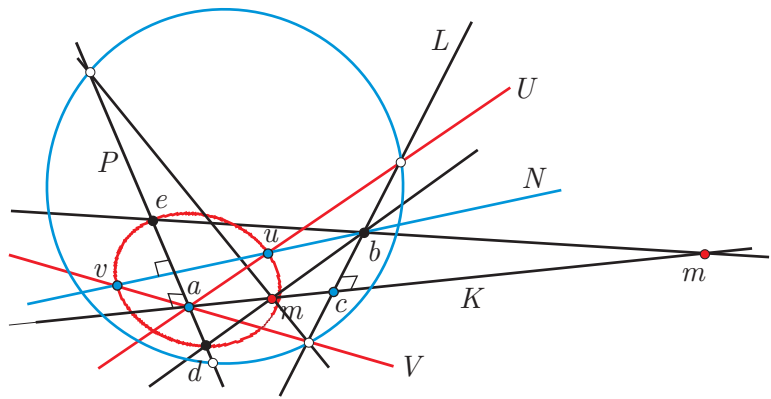}%
\caption{A variant on J. Bolyai's construction of the limiting lines from $a$
to $L$}%
\end{figure}
%EndExpansion

Start by constructing the altitude line $K$ from $a$ to $L,$ meeting $L$ at
$c,$ then the parallel line $P$ through $a$ to $L,$ namely that line
perpendicular to $K.$ Now let $m$ denote the midpoints of $\overline{ac},$
there are either two such points or none. If there are two, choose any point
$b$ on $L,$ construct the altitude $N$ to $P$ through $b,$ and reflect $b$ in
both midpoints $m$ to get $d$ and $e$ on $P.$ The side $\overline{ed}$ has $a$
a midpoint, and the hyperbolic circle centered at a through $d$ and $e$ meets
$N$ at the points $u$ and $v.$ Then $U\equiv au$ and $V\equiv av$ are the
required limiting lines as shown.

This construction seems to not be possible with only a straightedge, as we use
a hyperbolic circle; this would correspond to the fact that there are two
solutions. The question of what can and cannot be constructed with only a
straightedge seems also an interesting one.

\section{Canonical points}

Both the Canonical points theorem in this section and the Jumping Jack theorem
of the next section involve \textit{cubic relations} between certain
quadrances. I predict both will open up entirely new directions in hyperbolic geometry.

The Canonical points theorem has rather many aspects, one of which is a
classical theorem of projective geometry.

\begin{theorem}
[Canonical points]Suppose that $\alpha_{1}$ and $\alpha_{2}$ are distinct null
points, and that $x_{3}$ and $y_{3}$ are points lying on $\alpha_{1}\alpha
_{2}.$ For any third null point $\alpha_{3},$ and any point $b_{1}$ lying on
$\alpha_{2}\alpha_{3},$ define $x_{2}=\left(  \alpha_{1}\alpha_{3}\right)
\left(  y_{3}b_{1}\right)  $ and $y_{2}=\left(  \alpha_{1}\alpha_{3}\right)
\left(  x_{3}b_{1}\right)  $. Similarly for any point $b_{2}$ lying on
$\alpha_{1}\alpha_{3}$ define $x_{1}=\left(  \alpha_{2}\alpha_{3}\right)
\left(  y_{3}b_{2}\right)  $ and $y_{1}=\left(  \alpha_{2}\alpha_{3}\right)
\left(  x_{3}b_{2}\right)  $. Then $b_{3}\equiv\left(  x_{1}y_{2}\right)
\left(  x_{2}y_{2}\right)  $ lies on $\alpha_{1}\alpha_{2}.$ Now define points%
\[
c_{1}=\left(  x_{2}x_{3}\right)  \left(  y_{2}y_{3}\right)  \qquad
c_{2}=\left(  x_{1}x_{3}\right)  \left(  y_{1}y_{3}\right)  \qquad
c_{3}=\left(  x_{1}x_{2}\right)  \left(  y_{1}y_{2}\right)
\]
and corresponding points
\begin{align*}
z_{3}  &  =\left(  c_{1}b_{1}\right)  \left(  \alpha_{1}\alpha_{2}\right)
\qquad w_{2}=\left(  c_{1}b_{1}\right)  \left(  \alpha_{1}\alpha_{3}\right) \\
z_{1}  &  =\left(  c_{2}b_{2}\right)  \left(  \alpha_{2}\alpha_{3}\right)
\qquad w_{3}=\left(  c_{2}b_{2}\right)  \left(  \alpha_{1}\alpha_{2}\right) \\
z_{2}  &  =\left(  c_{3}b_{3}\right)  \left(  \alpha_{1}\alpha_{3}\right)
\qquad w_{1}=\left(  c_{3}b_{3}\right)  \left(  \alpha_{2}\alpha_{3}\right)  .
\end{align*}
Then $z_{3}$ and $w_{3}$ depend only on $x_{3}$ and $y_{3},$ and not on
$\alpha_{3},b_{1}$ and $b_{2}.$ Furthermore $b_{1},z_{2},w_{3}$ are collinear,
as are $b_{2},z_{3},w_{1},$ and $b_{3},z_{1},w_{2}$.%
%TCIMACRO{\FRAME{fhFU}{4.2983in}{3.3927in}{0pt}{\Qcb{Canonical points theorem:
%$\overline{x_{3}y_{3}}$ determines $\overline{z_{3}w_{3}}$}}%
%{\Qlb{TriangleNull Perspective}}{nulltriangleperpsective1.eps}%
%width 4.2983in;  height 3.3927in;  depth 0pt;  original-width 3.374in;
%original-height 2.659in;  cropleft "0";  croptop "1";  cropright "1";
%filename 'NullTrianglePerpsective1.eps';file-properties "XNPEU";}} }%
%BeginExpansion
\begin{figure}[h]%
\centering
\includegraphics[
height=3.3927in,
width=4.2983in
]%
{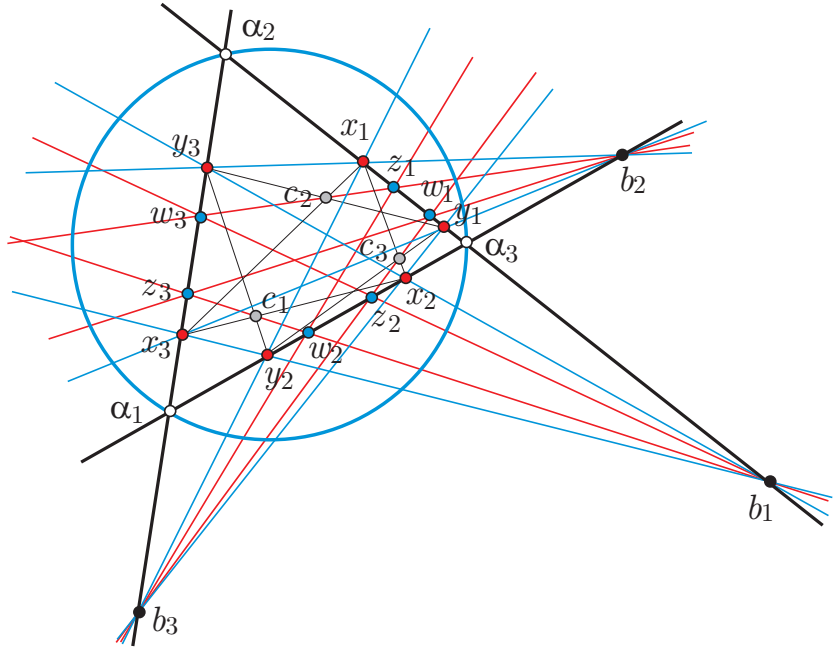}%
\caption{Canonical points theorem: $\overline{x_{3}y_{3}}$ determines
$\overline{z_{3}w_{3}}$}%
\label{TriangleNull Perspective}%
\end{figure}
%EndExpansion

\end{theorem}

In particular note that the theorem implies that any two points $x$ and $y$
whose join passes through two null points determine canonically two points $z
$ and $w$ lying on $xy$ in this fashion. We call $z$ and $w$ the
\textbf{canonical points} of $x$ and $y.$ In Figure
\ref{TriangleNull Perspective} $z_{3}$ and $w_{3}$ are the canonical points of
$x_{3}$ and $y_{3},$ while $z_{1}$ and $w_{1}$ are the canonical points of
$x_{1}$ and $y_{1},$ and $z_{2}$ and $w_{2}$ are the canonical points of
$x_{2}$ and $y_{2}$.

\begin{theorem}
[Canonical points cubic]With notation as above, the quadrances $q\equiv
q\left(  x_{3},y_{3}\right)  $ and $r\equiv q\left(  z_{3},w_{3}\right)  $
satisfy the cubic relation%
\begin{equation}
\left(  q-4r\right)  ^{2}=8qr\left(  2r-q\right)  . \label{Canonical}%
\end{equation}

\end{theorem}

\begin{figure}[ptb]%
\centering
\includegraphics[
height=2.4873in,
width=2.5247in
]%
{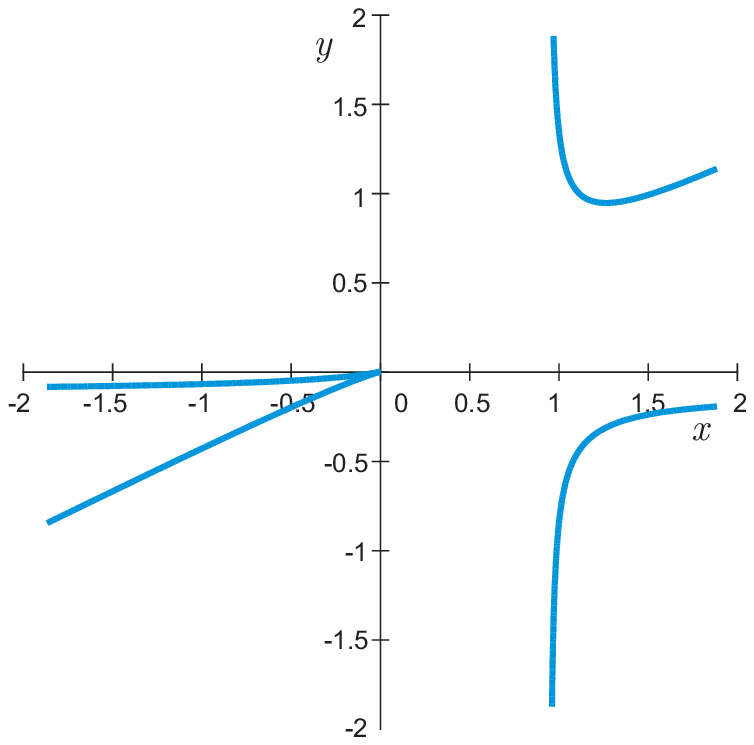}%
\caption{The Canonical points cubic: $\left(  x-4y\right)  ^{2}=8xy\left(
2y-x\right)  $}%
\label{CanonicalPointsCubic}%
\end{figure}

We call the algebraic curve%
\[
\left(  x-4y\right)  ^{2}=8xy\left(  2y-x\right)
\]
the \textbf{Canonical points cubic}. The graph is shown in Figure
\ref{CanonicalPointsCubic}. It is perhaps interesting that the point $\left[
9/8,9/8\right]  $ is the apex of one of the branches of this algebraic curve.%
%TCIMACRO{\FRAME{ftbpFU}{2.5247in}{2.4873in}{0pt}{\Qcb{The Canonical points
%cubic: $\left(  x-4y\right)  ^{2}=8xy\left(  2y-x\right)  $}}%
%{\Qlb{CanonicalPointsCubic}}{canonicalpointscubic.EPS}%
%width 2.5247in;  height 2.4873in;  depth 0pt;  original-width 2.8535in;
%original-height 2.8103in;  cropleft "0";  croptop "1";  cropright "1";
%cropbottom "0";  filename 'CanonicalPointsCubic.EPS';file-properties "XNPEU";}%
%} }%
%BeginExpansion
%EndExpansion

\section{The Jumping Jack theorem}

Here is my personal favourite theorem. Although one can give a computational
proof of it, the result begs for a conceptual framework that explains it, and
points to other similar facts (if they exist!)

\begin{theorem}
[Jumping Jack]\textit{Suppose that }$\overline{\alpha_{1}\alpha_{2}\alpha
_{3}\alpha_{4}}$\textit{\ is a quadrangle of null points, with }%
$g\equiv\left(  \alpha_{1}\alpha_{3}\right)  \left(  \alpha_{2}\alpha
_{4}\right)  $ \textit{a diagonal point, and let }$L$\textit{\ be any line
through }$g.$\textit{\ Then for an arbitrary null point }$\alpha_{5}%
,$\textit{\ define the meets }$x\equiv\left(  \alpha_{1}\alpha_{3}\right)
\left(  \alpha_{4}\alpha_{5}\right)  $,\textit{\ }$y\equiv L\left(  \alpha
_{4}\alpha_{5}\right)  $,\textit{\ }$z\equiv\left(  \alpha_{2}\alpha
_{4}\right)  \left(  \alpha_{3}\alpha_{5}\right)  $\textit{\ and }$w\equiv
L\left(  \alpha_{3}\alpha_{5}\right)  $\textit{. If }$r\equiv q\left(
x,y\right)  $\textit{\ and }$s\equiv q\left(  z,w\right)  $\textit{\ then}%
\[
16rs\left(  3-4\left(  s+r\right)  \right)  =1.
\]

\end{theorem}

%

%TCIMACRO{\FRAME{fhFU}{5.0959in}{2.0564in}{0pt}{\Qcb{Jumping Jack theorem:
%$16rs\left(  3-4\left(  s+r\right)  \right)  =1$}}{\Qlb{JumpingJack}%
%}{jumpingjackboth.EPS}{\special{ language "Scientific Word";  type "GRAPHIC";
%width 5.0959in;  height 2.0564in;  depth 0pt;  original-width 4.8385in;
%original-height 1.9369in;  cropleft "0";  croptop "1";  cropright "1";
%cropbottom "0";  filename 'JumpingJackBoth.EPS';file-properties "XNPEU";}} }%
%BeginExpansion
\begin{figure}[h]%
\centering
\includegraphics[
height=2.0564in,
width=5.0959in
]%
{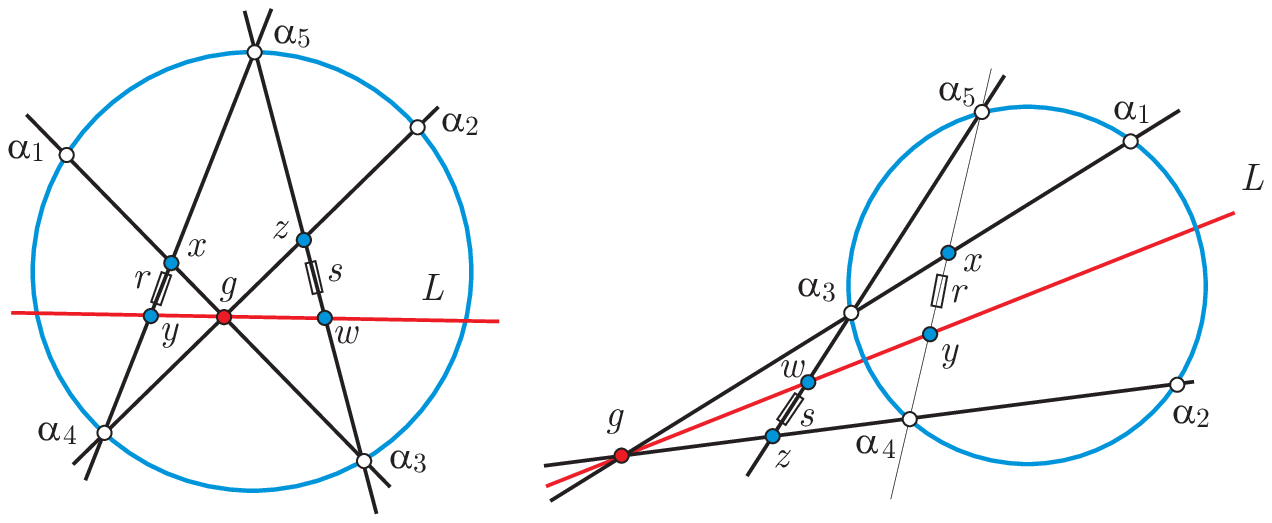}%
\caption{Jumping Jack theorem: $16rs\left(  3-4\left(  s+r\right)  \right)
=1$}%
\label{JumpingJack}%
\end{figure}
%EndExpansion

We call the algebraic curve
\[
16xy\left(  3-4\left(  x+y\right)  \right)  =1
\]
the \textbf{Jumping Jack cubic}. The Jumping Jack theorem shows that it has an
infinite number of rational solutions, which include a parametric description
with 6 independent parameters.

\begin{figure}[h]%
\centering
\includegraphics[
height=7.0116cm,
width=7.1171cm
]%
{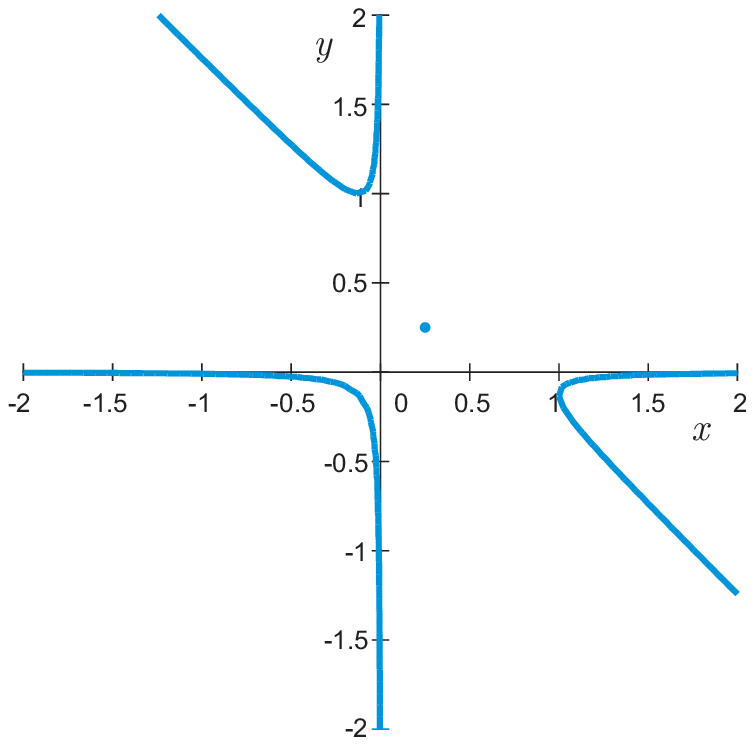}%
\caption{Jumping Jack cubic: $16xy\left(  3-4\left(  x+y\right)  \right)  =1$}%
\label{JumpingJackCubic1}%
\end{figure}

The graph is shown in Figure \ref{JumpingJackCubic1}. Note the isolated
solution $\left[  1/4,1/4\right]  $, which is the centroid of the trilateral
formed by the three asymptotes.%
%TCIMACRO{\FRAME{fhFU}{7.1171cm}{7.0116cm}{0pt}{\Qcb{Jumping Jack cubic:
%$16xy\left(  3-4\left(  x+y\right)  \right)  =1$}}{\Qlb{JumpingJackCubic1}%
%}{jumpingjackcubic.eps}{\special{ language "Scientific Word";
%type "GRAPHIC";  maintain-aspect-ratio TRUE;  display "USEDEF";
%valid_file "F";  width 7.1171cm;  height 7.0116cm;  depth 0pt;
%original-width 7.4566cm;  original-height 7.3438cm;  cropleft "0";
%filename 'JumpingJackCubic.eps';file-properties "XNPEU";}} }%
%BeginExpansion
%EndExpansion

\section{Conclusion}

Universal hyperbolic geometry provides a new framework for a classical
subject. It provides a more logical foundation for this geometry, as now
analysis is not used, but only high school algebra with polynomials and
rational functions. The main laws of trigonometry require only quadratic
equations for their solutions. Theorems extend now beyond the familiar
interior of the unit disk, and also to geometries over finite fields. Although
we have not stressed this, it turns out that almost all the theorems we have
described also hold in elliptic geometry! That is because the algebraic
treatment turns out to be essentially independent of the projective quadratic
form in the three dimensional space that is implicitly used to set up the
theory in (\ref{LineEquationDef}). We have shown how many classical results
can be enlarged to fit into this new framework, and also described new and
interesting results.

So there are many new opportunities for researchers to make essential
discoveries at this early stage of the subject. When it comes to hyperbolic
geometry, we are all beginners now.

\end{document}